\documentclass[10pt]{article}
\usepackage{amsmath,amssymb,latexsym,epsfig,subfigure,multirow,color}
\begin{document}

\title{Finite Volume Formulation of the MIB Method for Elliptic Interface Problems}

\author {Yin Cao$^{1}$, Bao Wang$^{1}$, Kelin Xia$^{1}$, and Guowei Wei$^{1,2}$
\footnote{
Corresponding author. Tel: (517)353 4689, Fax: (517)432 1562, Email: wei@math.msu.edu }\\
\\
\small \it     $^1$Department of Mathematics, Michigan State University, East Lansing, MI 48824, USA \\
\small \it     $^2$Department of Electrical and Computer Engineering, \\
\small \it          Michigan State University, East Lansing, MI 48824, USA }

\date{\today}

\maketitle

\begin{abstract}
The matched interface and boundary  (MIB) method has a proven ability for delivering the second order accuracy in handling elliptic interface problems with arbitrarily complex interface geometries. However, its collocation formulation requires relatively high solution regularity. Finite volume method (FVM) has its merit in dealing with conservation law problems and its integral formulation works well with relatively low solution regularity. We propose an MIB-FVM to take the advantages of both MIB and FVM for solving elliptic interface problems.  We construct the proposed method on  Cartesian meshes with  vertex-centered control volumes. A  large number of numerical experiments are designed to validate the present method in both two dimensional (2D) and three dimensional (3D) domains. It is found that the proposed MIB-FVM achieves the second order convergence for  elliptic interface problems with  complex interface geometries in both $L_{\infty}$ and $L_2$ norms.

\end{abstract}

\newpage
\section{Introduction}

Elliptic partial differential equations (PDEs) with discontinuous coefficients and singular source terms, commonly referred to as elliptic interface problems, occur in many applications, including fluid dynamics \cite{Fadlun:2000,Iaccarino:2003,Lee:2003a,Layton:2009}, material science \cite{Horikis:2006}, electromagnetics \cite{Hadley:2002,Hesthaven:2003,Kafafy:2005,Zhao:2004}, biological systems, \cite{Geng:2007a,Liu:2006,Yu:2007,Zhou:2008b,DuanChen:2011a} and heat or mass transfer \cite{Francois:2003}. Since the pioneer work of Peskin in 1977 \cite{Peskin:1977}, lots of attention has been paid to this field in the past few decades \cite{Ameur:2004,Berthelsen:2004,Dumett:2003,Fogelson:2001,Hou:1997,Hou:2005a,Jin:2002,Johansen:1998,Kandilarov:2003,Lee:2003a,Linnick:2005,Lombard:2003,Schulz:2002,Sethian:2001,Tornberg:2004,Vande-Voorde:2004,Morgenthal:2007,Wiegmann:2000,Cai:2003}.

An elliptic interface problem is formulated as  elliptic PDEs  defined  on piecewise-smooth  subdomains, which are coupled together via  interface conditions, such as given jumps in solution and flux across the domain interface. Without the interface, a large variety of convergent numerical schemes are available for solving the elliptic PDEs, such as, finite difference method (FDM), finite element method (FEM), finite volume method (FVM), wavelet, radial basis functions, meshless and spectral methods. However, for elliptic interface problems, the direct application of the aforementioned schemes cannot yield a  convergent solution.

 Peskin proposed the immersed boundary method (IBM) \cite{Griffith:2005,Lai:2000,Peskin:1977} in order to simulate flow pattern of blood in the heart, in which he approximated the singular sources on the interface. Since his work, many  numerical methods have been proposed and studied. In 1984, Mayo introduced a second order integral equation approach for Poisson's equation and biharmonic equation on irregular domains \cite{Mayo:1984,McKenney:1995}, in which the solution is extended to a rectangular region by using Fredholm integral equations. A fast Poisson solver is utilized to solve the resulting Fredholm integral equations on a rectangular region. Continuous derivatives have been assumed to evaluate the discrete Laplacian. This method can deal with jump conditions of $[u] \neq 0$ and $[u_n]=0$ when the Green's function is available. In 1994, a level set method in combination with the immersed boundary method \cite{MSussman:1994} is proposed by Osher and his coworkers in order to compute solutions to incompressible two-phase flow.  This level set method is easy to implement, although it is of low  convergent rate. In the same year, immersed interface method (IIM) was proposed by Levique and Li \cite{LeVeque:1994,Li:2001,Adams:2002}, in which the interface conditions are incorporated into the finite difference scheme near the interface to achieve second order accuracy based on a Taylor expansion in a local coordinate system. It is a second order interface scheme that does not smear the interface jumps and is used to solve elliptic and parabolic interface problems alike. The resulting linear system is sparse, but not symmetric or positive definite. Also, second order derivatives of the interface jumps are needed. The IIM is widely used in practice and its extensions can be seen in Refs. \cite{Li:2001,Adams:2002}. For instance, in 2001, Li and Ito \cite{Li:2001}   constructed an IIM scheme with resulting linear matrix being diagonally dominate and its symmetric part being negative definite.


In addition to the IIM, a large class of finite difference based numerical methods have also been proposed for solving elliptic interface problem. The crucial idea of these methods is to incorporate jump conditions into the difference stencils near the interface to maintain high order local truncation error using Taylor expansion. Numerical schemes based on finite element or finite volume are also  developed \cite{Overmann:2006,Li:2003,Hellrung:2012}.  Finite element schemes are usually constructed by modifying the finite element basis near the interface. A few methods on unfitted mesh have been introduced based on the discontinuous Galerkin method using interior penalty technique to handle jump and flux conditions \cite{Bastian:2009,Guyomarch:2009}.
The weak Galerkin finite element method has also been developed for elliptic interface problem  \cite{LMu:2013a}. Other interesting methods developed in recent decades include the ghost fluid method (GFM) by Fedkiw, Osher and coworkers \cite{Fedkiw:1999,Liu:2000}. The second order convergence of finite difference formalism of interface methods is proved in 2006 by Beale and Layton for smooth interface \cite{Beale:2006}, whereas convergence analysis of most   elliptic interface schemes is yet to be done.


In the past decade, we have dedicated ourselves to designing accurate and robust numerical schemes for solving   elliptic interface problems. In 2006, the matched interface and boundary (MIB) method \cite{Yu:2007c,Yu:2007a,Zhou:2006c,Zhou:2006d} was proposed, motivated by many practical needs, such as  optical molecular imaging \cite{DuanChen:2010b}, nano-electronic devices, \cite{DuanChen:2010b},  vibration analysis of plates \cite{SNYu:2009}, wave propagation \cite{SZhao:2010a,SZhao:2008a},   geodynamics   \cite{YCZhou:2012a} and electrostatic potential in proteins \cite{Zhou:2008b,Yu:2007,Geng:2007a, DuanChen:2011a}. One important feature of the MIB method is its extension of computational domains by using the so called fictitious values, an strategy developed in our earlier methods for handling boundaries \cite{JCPWei:1999,Wei:2002f} and interfaces \cite{Zhao:2004}. As such, standard  central finite difference schemes can be employed to discretize differential operators as if there were no interface. Another unique feature of the MIB method is to repeatedly enforce only the lowest order jump conditions to achieve higher order convergence, which is of critical importance for the robustness of the method to deal with arbitrarily complex interface geometries. Higher-order jump conditions must involve higher order derivatives and/or cross derivatives. Therefore, to approximate higher order derivatives and/or cross derivatives, one must utilize larger stencils, which is unstable for constructing high-order interface schemes and cumbersome for complex interface geometries.
The other distinct feature of the MIB method is its dimensional splitting. To enforce a  2D or 3D interface jump condition, we divide the problem into multi 1D ones and resolve a 1D problem at a time, if possible. This approach enables us to come up with a systematic procedure to resolve high dimensional interface problems.
Finally,  based on high order Lagrange polynomials, the MIB method is of arbitrarily high order in principle. For example, MIB schemes up to 16th order accurate have been constructed for simple interface geometries 1D and 2D domains \cite{Zhao:2004,Zhou:2006c}, and sixth-order accurate MIB schemes have been developed for complex interfaces in 2D \cite{Zhou:2006c} and 3D domains \cite{Yu:2007a,Zhou:2006c}. Recently we have constructed an adaptively  deformed  mesh  based MIB \cite{KLXia:2012a}  and a Galerkin formulation of MIB \cite{KLXia:2014e,KLXia:2014f} to improve MIB's capability of solving realistic problems. A comparison of the GFM, IIM and MIB methods can be found in in Refs. \cite{Zhou:2006c,Zhou:2006d}.

The MIB method has been used in our earlier works for solving many scientific and engineering problems, such as Poisson-Boltzmann equation (PBE) \cite{Zhou:2008b,Yu:2007,Geng:2007a, DuanChen:2011a} for describing the electrostatic potential in proteins. To our best knowledge, the MIB method is the only method that has demonstrated the second order accuracy for solving Poisson or PB equation with realistic protein surfaces with geometric singularities \cite{Yu:2007,Geng:2007a, DuanChen:2011a} and for solving Poisson equation with multiple material interfaces \cite{KLXia:2011}. It can also be used to solve the Helmholtz equation for wave scattering and propagation in inhomogeneous media. A fourth order MIB scheme for the Helmholtz equation with arbitrarily curved interfaces has been proposed by Zhao \cite{SZhao:2010a,SZhao:2008a}. Another example is the geodynamics where the Navier-Stokes equations with discontinuous viscosity and density is to be solved.  Zhou {\it et al} have  developed a second order accurate MIB method to solve this problem on non-staggered Cartesian grids \cite{YCZhou:2012a}.  Furthermore, the elastic interface problem in both 2D and 3D domains with arbitrarily complex interface geometry was also addressed by the MIB method \cite{BaoWang:2015c,BaoWang:2015d}. In the past few years, the MIB method has also been applied to the optical molecular imaging \cite{DuanChen:2010}, nano-electronic devices, \cite{DuanChen:2010b} and vibration analysis of plates \cite{SNYu:2009}.

Due to the importance and complexity of interface problems, it is still urgent to develop  methods that are more accurate, robust and require less computational cost. Finite volume method (FVM) is known for its ability to better conserve mass and flux in conservation law problems. This property is fundamental for the simulation of many physical models, e.g.,  oil recovery simulation and computational fluid dynamics in general. Additionally, compared with  finite difference method which takes a collocation formulation, FVM is able to deal with solutions with relatively lower regularity. In the present work, we propose the finite volume formulation of the MIB method (MIB-FVM). The motivation for the proposed MIB-FVM is to inherit the merits of both methods, namely, the capability of the MIB method in handling complex interface geometries and the advantage of FVM for conservation laws and low solution regularity.

The rest of this paper is organized as the follows. In Section \ref{theory_and_alg}, the general theory of the MIB based finite volume formulation is briefly discussed.  The theoretical formulation and the computational algorithms are given as well. In Section \ref{mib}, we present the strategy for treating complex interface geometries. In Section \ref{Numericalstudies}, the proposed MIB-FVM is validated by benchmark tests, such as  6-petal flower and jigsaw puzzle-like shape in two dimensional (2D) space, and   sphere, ellipsoid, cylinder, flower-based cylinder and torus in three dimensional (3D) space. Solutions of less regularity ($H^2$ continuous) are also tested for the spherical and ellipsoidal interface in 3D. This paper ends  with a conclusion.

\section{Matched interface and boundary - finite volume method (MIB-FVM)} \label{theory_and_alg}

\begin{figure}
\begin{center}
\includegraphics[width=0.5\textwidth]{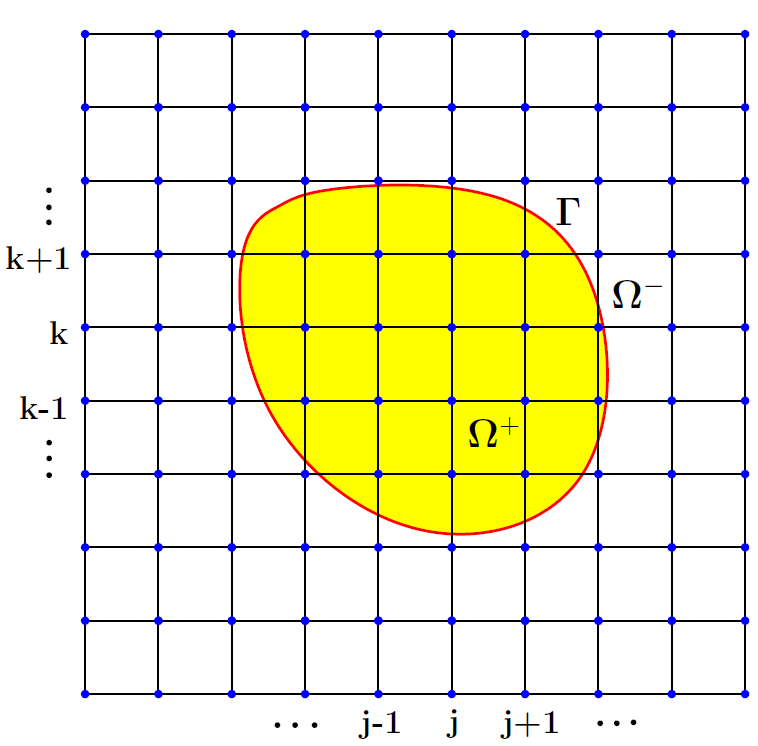}
\end{center}
\caption{Illustration of a 3D interface problem at a cross section. The whole domain $\Omega$ is divided by  interface $\Gamma$ into two subdomains, $\Omega^+$ and $\Omega^-$. The mesh is represented by solid lines. The grid points are represented by blue dots.}
\label{interface}
\end{figure}

Let us consider an open bounded domain $\Omega \subset {\mathbb R}^3$ with a given interface $\Gamma$, which separates the domain into two subdomains,  $\Omega^+$ and  $\Omega^-$ as illustrated in Figure \ref{interface}.

The boundary $\partial\Omega$ and interfaces $\Gamma$ may be non-smooth. The interface can be characterized by a piecewise smooth level-set function $\varphi\in \Omega$, such that $\Gamma=\{{\bf x}|\varphi({\bf x}) =0, \forall{\bf x}\in \Omega\}$.  As such, two subdomains can be given by $\Omega^+=\{{\bf x}|\varphi({\bf x})  < 0, \forall{\bf x}\in \Omega\}$ and $\Omega^-=\{{\bf x}|\varphi({\bf x}) > 0, \forall{\bf x}\in \Omega\}$. The elliptic interface problem can be formulated as
\begin{eqnarray}\label{theoryeq1}
-	\nabla\cdot \beta({\bf x})\nabla u({\bf x}) &=& g({\bf x}), ~\forall {\bf x}\in \Omega \\
                                                   u&=&g_b, ~\forall {\bf x}~  {\rm on} ~ \partial \Omega,
\end{eqnarray}
where $g({\bf x})$ is a piecewise continuous function, $g_b$ is the boundary value, and  $\beta({\bf x})$ is a  variable coefficient that is discontinuous across the interface $\Gamma$.  As a result, two jump conditions are required to make the problem well posed
\begin{eqnarray}
\label{jumpcondition1}
[u]        &=  u^+ - u^-     = \Phi,~     \forall {\bf x} ~{\rm on} ~\Gamma
\end{eqnarray}
and
\begin{eqnarray}
\label{jumpcondition2}
[\beta u_n]    &=   \beta^+ u_n^+ - \beta^- u_n^-   =  \Psi,~  \forall {\bf x} ~{\rm on} ~\Gamma,
\end{eqnarray}
where $u^+, u^+_n$ and $\beta^+$ denote their limiting value from the $\Omega^+$ side of the interface $\Gamma$, and   $u^-, u^-_n$ and $\beta^-$ denote their limiting value from the $\Omega^-$ side of the interface $\Gamma$. The derivatives $u^+_n$ and $u^-_n$ are evaluated along the outer normal direction on the interface. Here  $\Phi({\bf x})$ and $\Psi({\bf x})$ is at least $C^1$ continuous. Equations (\ref{theoryeq1})-(\ref{jumpcondition2}) define the elliptic interface problem to be solved in the present work.

\subsection{Introduction to the MIB finite volume method}

\begin{figure}
\begin{center}
\begin{tabular}{cc}
\includegraphics[width=0.5\textwidth]{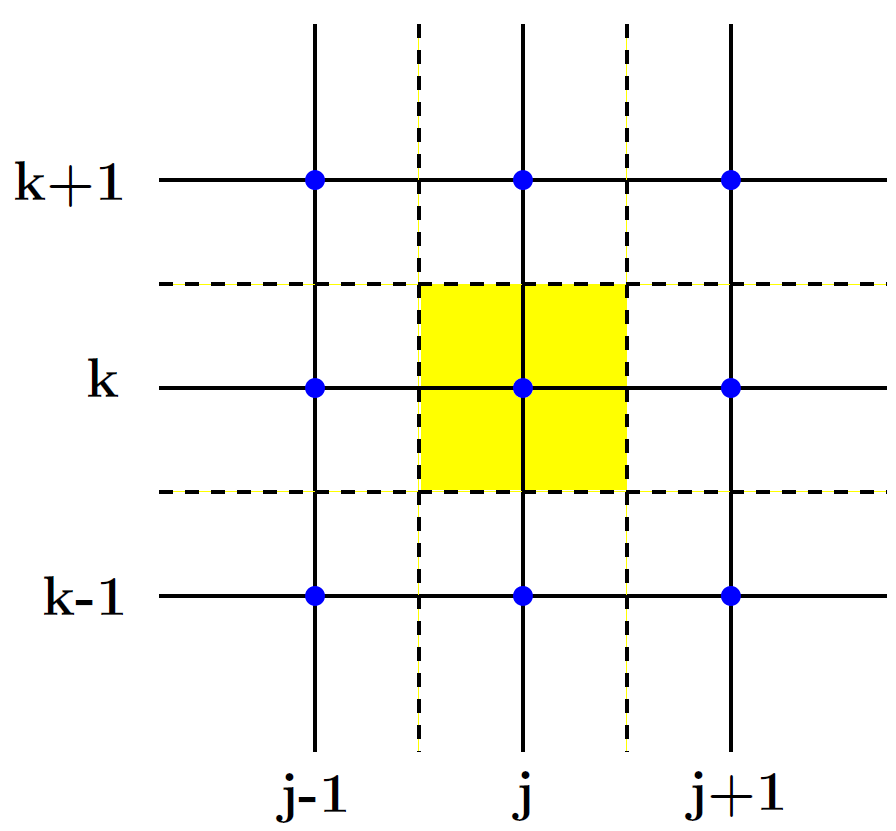} &
\includegraphics[width=0.5\textwidth]{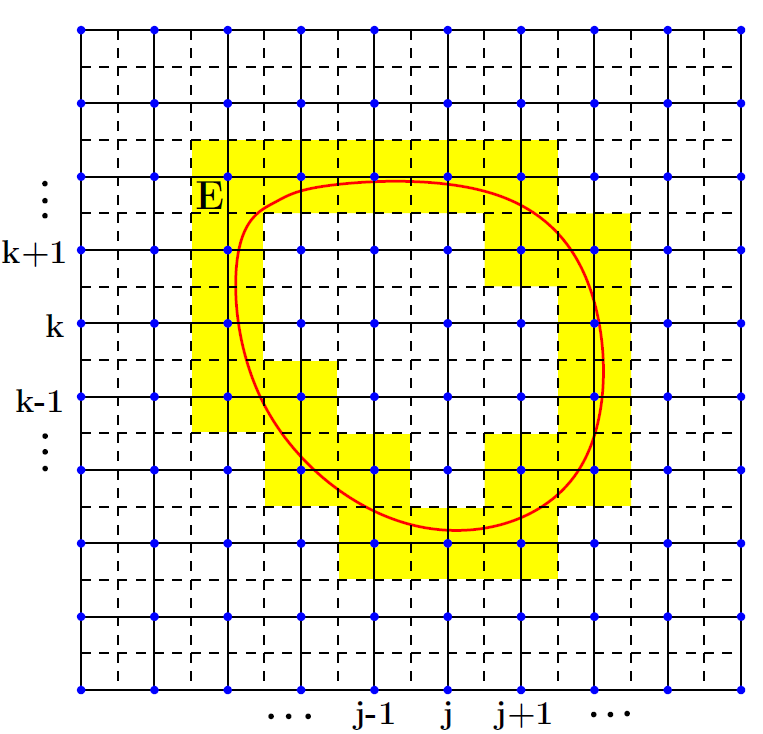}
\end{tabular}
\end{center}
\caption{Left chart: Illustration of the 3D control volume at a cross section that is perpendicular to the $x$-axis and intersect with the $x$-axis at the $i$th grid point. The yellow region denotes the control volume associated with the grid point $(i,j,k)$. The grid points are denoted by blue dots. The mesh is represented by solid lines and the dual mesh is represented by the dashed lines. Right chart: Illustration of the 3D irregular domain $E$ at a cross section. The irregular domain is denoted by the yellow region. The grid points are denoted by blue dots. The mesh is represented by solid lines and the dual mesh is represented by dashed lines.}
\label{cv}
\end{figure}

To avoid the time-consuming mesh generation procedure, like the original MIB method, we employ the Cartesian mesh. We also employ the vertex-centered finite volume method, which means that we associates control volumes and unknowns to vertices. The control volumes are cubes whose centers are located at the grid points and whose sides intercept at the midpoints between grid points, which is illustrated by the left chart of Figure \ref{cv}.

Like the original MIB method, the grid points and and control volumes are classified into regular and irregular types. A grid point is said to be irregular if the standard central finite difference (CFD) scheme at the grid point near the interface requires grid point(s) from the other side of the interface. In Figure \ref{interface}, the interface $\Gamma$ and the different domains $\Omega$, $\Omega^+$ and $\Omega^-$ are depicted. As is shown in the figure, the complex interface inevitably cut through certain control volumes. A control volume is defined to be irregular if it is cut through by the interface, or stated differently, the vertices of the control volume locate on both side of the interface. The irregular domain $E$ is defined to be the union of all irregular control volumes, as is shown in the right chart of Figure \ref{cv}.

\begin{figure}
\begin{center}
\begin{tabular}{cc}
\includegraphics[width=0.5\textwidth]{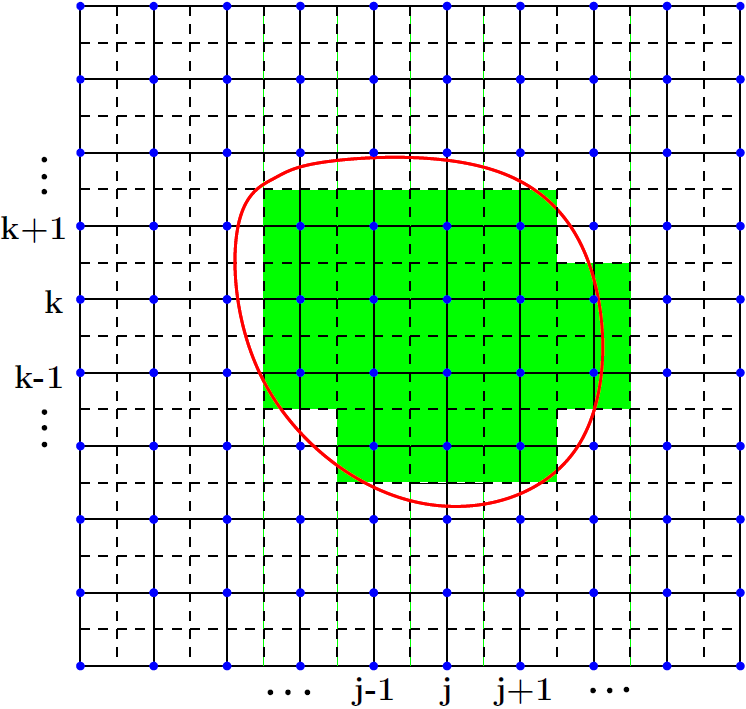} &
\includegraphics[width=0.5\textwidth]{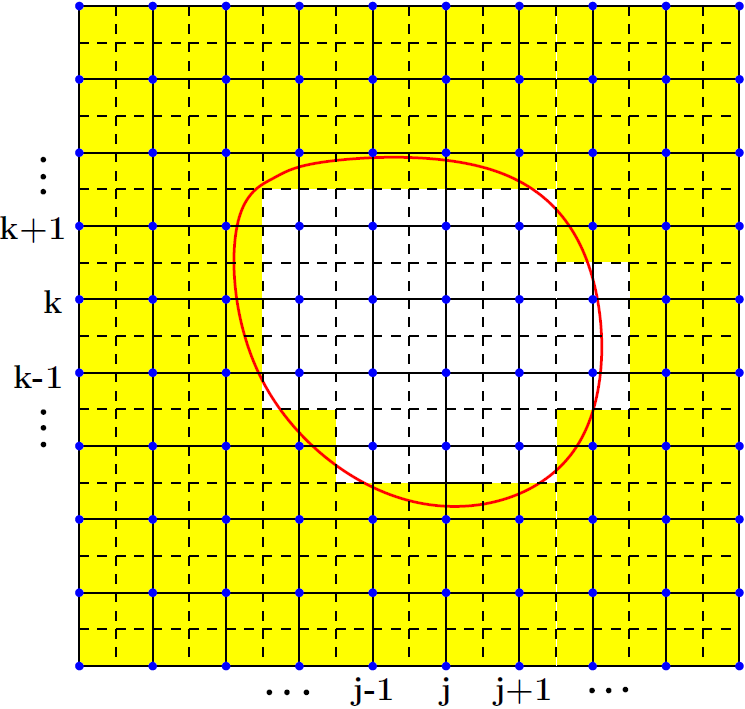}
\end{tabular}
\end{center}
\caption{Illustration of the domain re-division at a cross section. Left chart: domain $\Omega_o^+$ (in green), which is the union of all the positive control volumes; Right chart: domain $\Omega_o^-$ (in yellow), which is the union of all the negative control volumes. }
\label{domainredivision}
\end{figure}

In the present MIB-FVM, we further define each control volume to be a positive or negative control volume according to where its center is, i.e., a control volume is said a positive one if its center grid point locates in the positive region $\Omega^+$; similarly, a control volume is said a negative one if its center grid point locates in the negative region $\Omega^-$. We also define $\Omega_o^+$ and $\Omega_o^-$ to be the union of all positive and negative control volumes, respectively, as is shown in Figure \ref{domainredivision}. Domain $\Omega_o^+$ and $\Omega_o^-$ can be viewed as the generalized $\Omega^+$ and $\Omega^-$. First, domains $\Omega_o^+$ and $\Omega_o^-$ are disjoint by such definition; secondly, it can be easily seen that there is no inclusion relation between the origin subdomains $\Omega^+$ or $\Omega^-$ and the generalized subdomains $\Omega_o^+$ or $\Omega_o^-$.

To properly define the elliptic interface problem, we also need to generalize the discontinuous coefficient $\beta (\bf x)$ and the source term $ g(\bf x)$ over the generalized domains. To this end,  we define

\begin{align} \label{exbetapos}
\beta_{o}^+({\bf x})=
\begin{cases}
       \beta^+({\bf x}),                   \qquad & \forall {\bf x} \subset \Omega_{o}^+ \cap \Omega^+ \\
       \beta_{ex}^+({\bf x}),             \qquad & \forall {\bf x} \subset \Omega_{o}^+ \cap {\Omega^+}^{c}
\end{cases}
\end{align}
and
\begin{align} \label{exbetaneg}
\beta_{o}^-({\bf x})=
\begin{cases}
       \beta^-({\bf x}),                   \qquad & \forall {\bf x} \subset \Omega_{o}^- \cap \Omega^- \\
       \beta_{ex}^-({\bf x}),             \qquad & \forall {\bf x} \subset \Omega_{o}^- \cap {\Omega^-}^{c},
\end{cases}
\end{align}
where $\beta_{ex}^+({\bf x})$ and $\beta_{ex}^-({\bf x})$ are the smooth extensions of the coefficients $\beta^+({\bf x})$ and $\beta^-({\bf x})$ over domains $\Omega_{o}^+$ and $\Omega_{o}^-$, respectively. The superscript $c$ means complement.

Similarly, we define
\begin{align}
g_{o}^+({\bf x})=
\begin{cases}
       g^+({\bf x}),                   \qquad & \forall {\bf x} \subset \Omega_{o}^+ \cap \Omega^+ \\
       g_{ex}^+({\bf x}),             \qquad & \forall {\bf x} \subset \Omega_{o}^+ \cap {\Omega^+}^{c}
\end{cases}
\end{align}
and
\begin{align}
g_{o}^-({\bf x})=
\begin{cases}
       g^-({\bf x}),                   \qquad & \forall {\bf x} \subset \Omega_{o}^- \cap \Omega^- \\
       g_{ex}^-({\bf x}),             \qquad & \forall {\bf x} \subset \Omega_{o}^- \cap {\Omega^-}^{c},
\end{cases}
\end{align}
where $g_{ex}^+$ and $g_{ex}^-$ are the smooth extension of the functions $g^+({\bf x})$ and $g^-({\bf x})$ over domains $\Omega_{o}^+$ and $\Omega_{o}^-$, respectively. The superscript $c$ denotes complement.

By doing the generalization above, we can define the original interface problems on irregular control volumes. As the solution and its gradient can be discontinuous across the interface, to solve the interface problem, one needs to rigorously enforce the interface jump conditions. These conditions, interpreted by a set of equations, are utilized to determine fictitious values on generalized domains numerically. The detailed procedure of computing fictitious values is described in Section \ref{mib}. To achieve high order finite difference type of MIB schemes, the set of interface conditions are iteratively implemented to determine as many fictitious values as needed. This approach can also be used to create a high-order MIB-FVM.

\subsection{MIB finite volume formulation}
The MIB-FVM is based on the {\bf integral form} rather than the {\bf differential form} of Eq.  (\ref{theoryeq1}). The integral form of the equation can be obtained by taking integral on both sides of the equation over any control volume $V_i$ and using the divergence theorem on the left hand side
\begin{align} \label{integralform1}
\int_{\partial V_i}{\beta({\bf x})\nabla u({\bf x})\cdot{\bf n}}dS = \int_{V_i}{g({\bf x})}d{\bf x},
\end{align}
where ${\partial V_i}$ denotes the surface of the control volume $V_i$, and {\bf n} denotes the unit outer normal vector of ${\partial V_i}$.

The above equation can be written as
\begin{align} \label{integralform2}
\int_{\partial V_i} (\beta u_x, \beta u_y, \beta u_z) \cdot (n_1, n_2, n_3)dS = \int_{V_i}{g({\bf x})}d{\bf x},
\end{align}
where $(n_1, n_2, n_3)$ is the coordinate of the unit outer normal $\bf n$.

\begin{figure}
\begin{center}
\begin{tabular}{ccc}
\includegraphics[width=0.33\textwidth]{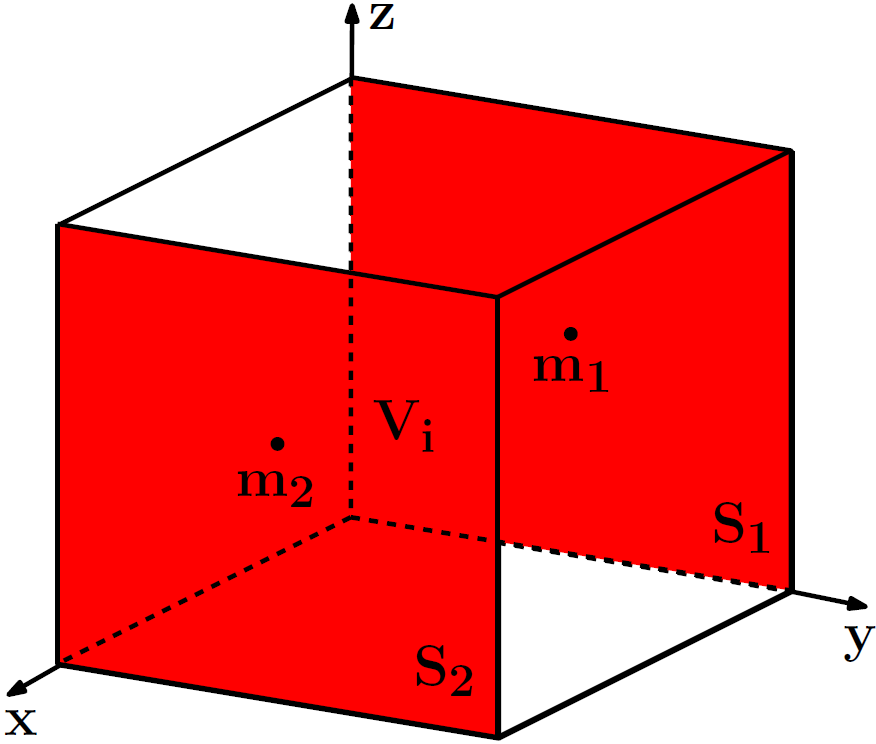} &
\includegraphics[width=0.33\textwidth]{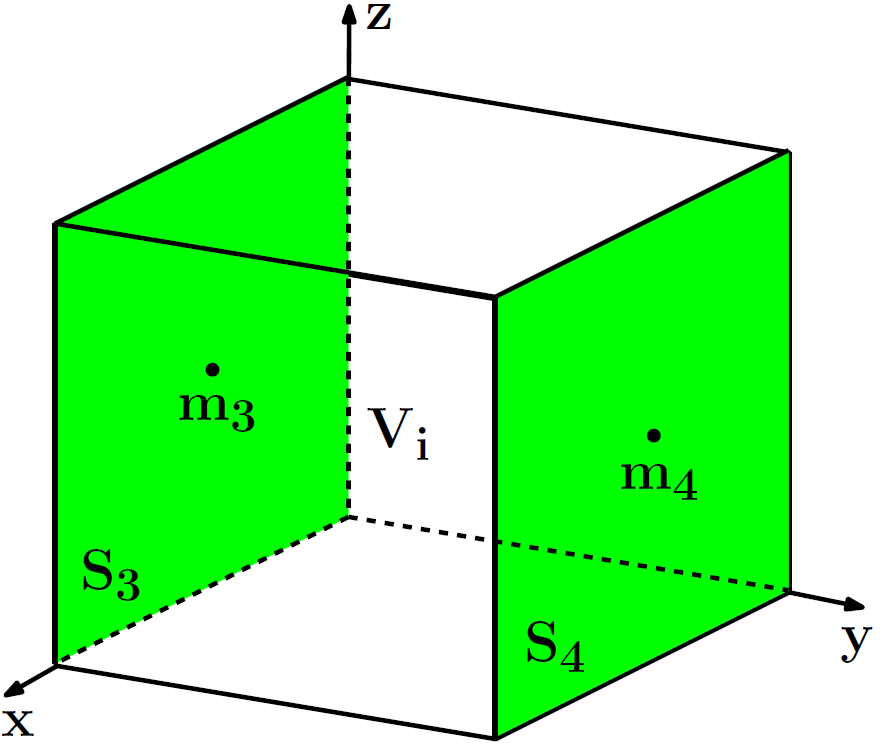} &
\includegraphics[width=0.33\textwidth]{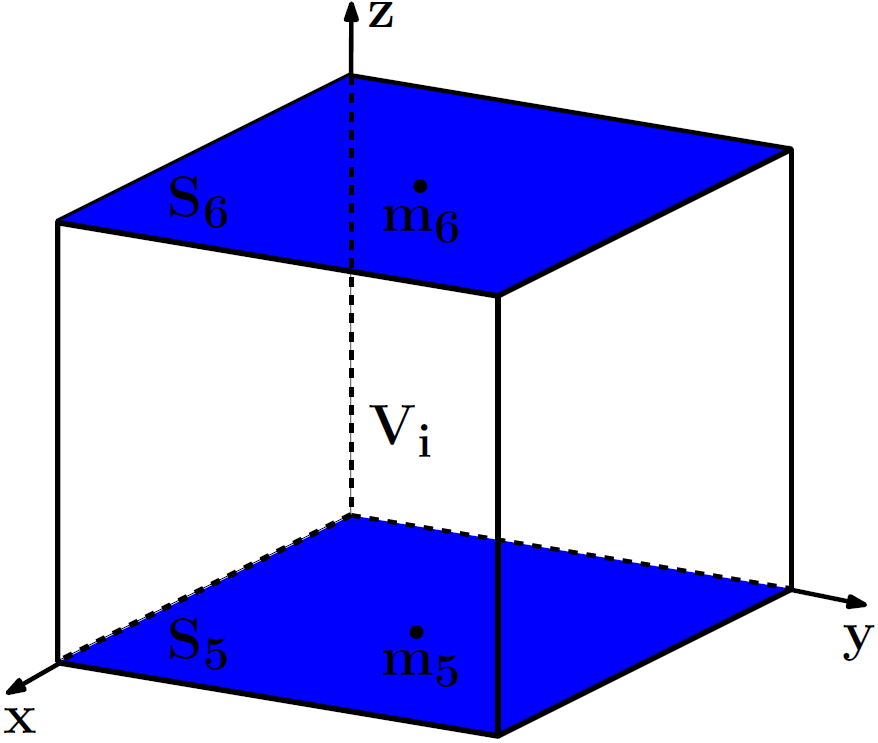}
\end{tabular}
\end{center}
\caption{Illustration of the faces of control volume $V_i$. Left chart: $S_1$ and $S_2$ are two faces along the $x$-direction; middle chart: $S_3$ and $S_4$ are two faces along the $y$-direction; right chart: $S_5$ and $S_6$ are two faces along the $z$-direction. $m_i$ is the center of $S_i$, where $i=1,2,...,6$.}
\label{faces}
\end{figure}

As is illustrated in Figure \ref{faces}, since the control volume $V_i$ is a cube, $\partial V_i$ consists of 6 faces, denoted by $\partial S_i$ where $i=1,2,...6$. For instance, $\partial S_1$ denotes the left face, whose unit outer normal is $(-1,0,0)$; $\partial S_2$ denotes the right face, whose unit outer normal is $(1,0,0)$, and so on.

As a result, Eq.  (\ref{integralform2}) can be further written as
\begin{align} \label{6surfaceintegral}
-\int_{\partial S_1} \beta u_x dS + \int_{\partial S_2} \beta u_x dS - \int_{\partial S_3} \beta u_y dS + \int_{\partial S_4} \beta u_y dS - \int_{\partial S_5} \beta u_z dS + \int_{\partial S_6} \beta u_z dS = \int_{V_i}{g({\bf x})}d{\bf x}.
\end{align}

\begin{figure}
\begin{center}
\begin{tabular}{cc}
\includegraphics[width=0.5\textwidth]{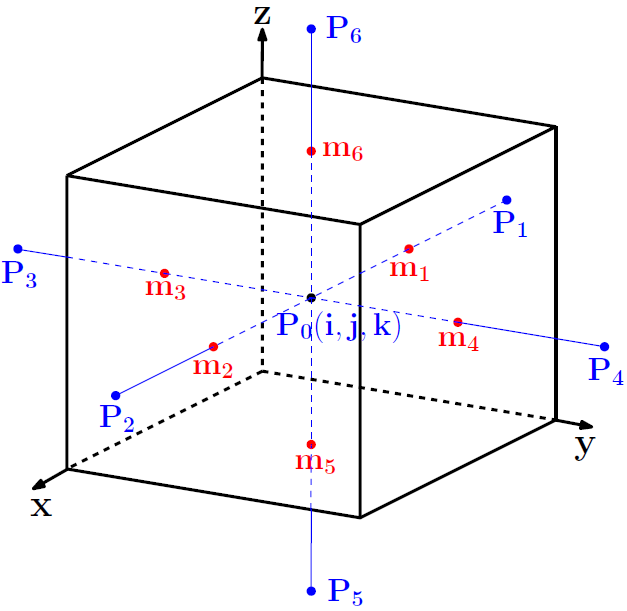} &
\includegraphics[width=0.5\textwidth]{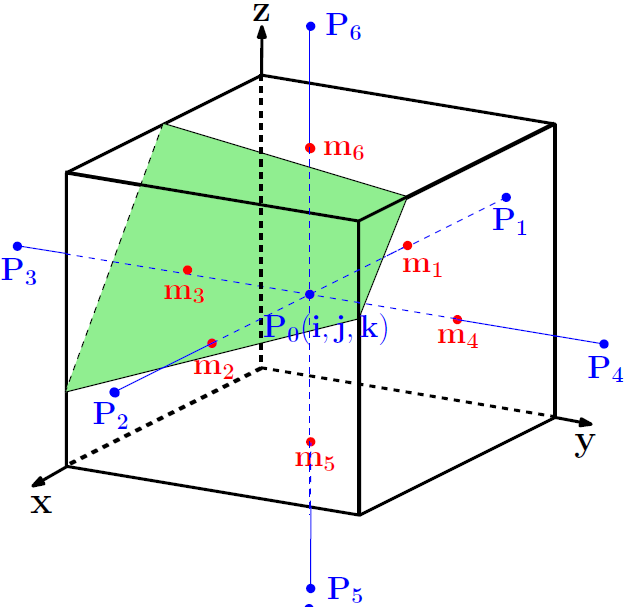}
\end{tabular}
\end{center}
\caption{Exemplary configurations for a regular control volume on the left and an irregular one on the right. Here $P_i, i=0,1,2,...,6$ in blue color denote the grid points, with $P_0$ being the center of the control volume. Here $m_i, i=1,2,...,6$ in red color denote the centers of the corresponding faces. For the right chart, the green area represents the 2D interface in 3D that divides the control volume into two parts.}
\label{controlvolume}
\end{figure}

The above surface and volume integrals are further approximated by Gauss quadrature rules. In practice, we use the   second order midpoint rule, thus the surface integrals on the left hand side are approximated by the function values at the face center, times the face area; the volume integral on the right hand side is approximated by the function value at the center, times the volume of the control volume. Based on Figure \ref{controlvolume} and the midpoint rule, Eq.  (\ref{6surfaceintegral}) can be further written as
\begin{align}
-(\beta u_x)_{m_1}  dydz + (\beta u_x)_{m_2}dydz - (\beta u_y)_{m_3}dxdz + (\beta u_y)_{m_4}dxdz - (\beta u_z)_{m_5}dxdy + (\beta u_z)_{m_6}dxdy = g_i dxdydz,
\end{align}
where $(\beta u_t)_{m_i}, t=x,y,z; i=1,2,...,6$  denote the function value of $(\beta u_t)$ at the $i$th midpoint $m_i$; $g_i$ denote the function value of $g({\bf x})$ at the center of the $i$th control volume $V_i$; and $dx, dy$ and $dz$ denote the size of partition in $x, y $ and $z$ direction, respectively.

For a regular control volume illustrated by the left chart of Figure \ref{controlvolume}, the above partial derivatives at the midpoints can be evaluated by a 3-point interpolation scheme, namely
\begin{align}
\begin{split}
(\beta u_x)_{m_1} \approx (\beta)_{m_1} \cdot (w_{1}^{1},w_{1}^{2},w_{1}^{3}) \cdot (u_{i-1,j,k},u_{i,j,k},u_{i+1,j,k})^T \\
(\beta u_x)_{m_2} \approx (\beta)_{m_2} \cdot (w_{2}^{1},w_{2}^{2},w_{2}^{3}) \cdot (u_{i-1,j,k},u_{i,j,k},u_{i+1,j,k})^T \\
(\beta u_y)_{m_3} \approx (\beta)_{m_3} \cdot (w_{3}^{1},w_{3}^{2},w_{3}^{3}) \cdot (u_{i,j-1,k},u_{i,j,k},u_{i,j+1,k})^T \\
(\beta u_y)_{m_4} \approx (\beta)_{m_4} \cdot (w_{4}^{1},w_{4}^{2},w_{4}^{3}) \cdot (u_{i,j-1,k},u_{i,j,k},u_{i,j+1,k})^T \\
(\beta u_z)_{m_5} \approx (\beta)_{m_5} \cdot (w_{5}^{1},w_{5}^{2},w_{5}^{3}) \cdot (u_{i,j,k-1},u_{i,j,k},u_{i,j,k+1})^T \\
(\beta u_z)_{m_6} \approx (\beta)_{m_6} \cdot (w_{6}^{1},w_{6}^{2},w_{6}^{3}) \cdot (u_{i,j,k-1},u_{i,j,k},u_{i,j,k+1})^T
\end{split}
\end{align}
where $w_i^j, i=1,2...,6; j=1,2,3$ are the interpolation weights, and $u_{i,j,k}$ denotes the function value at grid point  $(i,j,k)$.

For an irregular control volume illustrated by the right chart of Figure \ref{controlvolume}, to maintain the same convergence rate and accuracy, the above interpolations must be modified as
\begin{align}
\begin{split}
(\beta u_x)_{m_1} \approx (\beta_{o})_{m_1} \cdot (w_{1}^{1},w_{1}^{2},w_{1}^{3}) \cdot ({\tilde u}_{i-1,j,k},u_{i,j,k},{\tilde u}_{i+1,j,k})^T \\
(\beta u_x)_{m_2} \approx (\beta_{o})_{m_2} \cdot (w_{2}^{1},w_{2}^{2},w_{2}^{3}) \cdot ({\tilde u}_{i-1,j,k},u_{i,j,k},{\tilde u}_{i+1,j,k})^T \\
(\beta u_y)_{m_3} \approx (\beta_{o})_{m_3} \cdot (w_{3}^{1},w_{3}^{2},w_{3}^{3}) \cdot ({\tilde u}_{i,j-1,k},u_{i,j,k},{\tilde u}_{i,j+1,k})^T \\
(\beta u_y)_{m_4} \approx (\beta_{o})_{m_4} \cdot (w_{4}^{1},w_{4}^{2},w_{4}^{3}) \cdot ({\tilde u}_{i,j-1,k},u_{i,j,k},{\tilde u}_{i,j+1,k})^T \\
(\beta u_z)_{m_5} \approx (\beta_{o})_{m_5} \cdot (w_{5}^{1},w_{5}^{2},w_{5}^{3}) \cdot ({\tilde u}_{i,j,k-1},u_{i,j,k},{\tilde u}_{i,j,k+1})^T \\
(\beta u_z)_{m_6} \approx (\beta_{o})_{m_6} \cdot (w_{6}^{1},w_{6}^{2},w_{6}^{3}) \cdot ({\tilde u}_{i,j,k-1},u_{i,j,k},{\tilde u}_{i,j,k+1})^T
\end{split}
\end{align}
where $\beta_{o}$ denotes the generalization of $\beta$ defined in Eqs.   (\ref{exbetapos}) and (\ref{exbetaneg}), and ${\tilde{u}}$ is defined as
	\begin{equation*}
		{\tilde{u}}_{x,y,z}=
		\begin{cases}
		u_{x,y,z} & \mbox{if $(x,y,z)$ is on the same side of the interface with $(i,j,k)$ }\\
		f_{x,y,z} & \mbox{if $(x,y,z)$ is on different side of the interface with $(i,j,k)$}.
		\end{cases}
	\end{equation*}
	Here $f_{x,y,z}$ represents the fictitious value at grid point $(x,y,z)$, which is computed by MIB method. The detailed procedure of determining these fictitious values 		are described in the next section.

\section{Algorithms for determining fictitious values} \label{mib}
\subsection{Simplification of interface jump conditions}

To solve the above linear equation system, we need to determine the involved fictitious values first. This issue is discussed in the present section. We first describe the interface jump conditions, which are needed at each intersecting point of control volume edges and the interface.

As the interface normal direction varies along the interface, which is very troublesome from the computational point of view. We introduce a local coordinate system at each intersection point of the Cartesian mesh and the interface to treat different interface geometries systematically.  At a specific intersection point, the local coordinate system is denoted by $(\xi, \eta, \zeta)$, where $\xi$ denotes the normal direction and $\eta$ is on the $x-y$ plane. The relation of the local coordinates and the Cartesian coordinates is described in the following
\begin{align} \label{transformrelation}
    \left(\begin{array}{c}
           \xi\\
           \eta	\\
           \zeta\\
	      \end{array}\right)&={\bf P } \cdot \left(
	      \begin{array}{c}
           x\\
           y\\
           z\\
	      \end{array}\right)
\end{align}
where ${\bf P}$ is the transformation matrix
\begin{align} \label{transformationmatrix}
     {\bf P}&=\left(
	      \begin{array}{ccc}
           \sin\phi\cos\theta  &\sin\phi\sin\theta   &\cos\phi \\
           -\sin\theta           &\cos\theta             &0             \\
           -\cos\phi\cos\theta  &-\cos\phi\sin\theta   &\sin\phi  \\
	      \end{array}\right),
\end{align}
where $\theta$ and $\phi$ are the azimuth and zenith angles with respect to the normal direction $\bf n$, respectively.

Without generating higher order jump conditions, usually we differentiate Eq.  (\ref{jumpcondition1}) in two tangential directions in order to obtain two additional jump conditions $[u_{\eta}]=u_{\eta}^+ -u_{\eta}^-$, and  $[u_{\zeta}]=u_{\zeta}^+ -u_{\zeta}^-$. Thus, in the new coordinate system, the jump conditions in two tangential directions can be discretized as
\begin{align}
\label{3dbasiccondition3}&[u_{\eta}]=(-\sin\theta u_x^+  + \cos \theta u_y^+ ) - (-\sin \theta u_x^-  + \cos \theta u_y^- )  \\
\label{3dbasiccondition4}&[u_{\zeta}]=(-\cos\phi \cos\theta u_x^+ - \cos \phi \sin \theta u_y^+ +\sin \phi u_z^+ ) - (-\cos \phi \cos \theta u_x^-  -  \cos \phi \sin \theta u_y^- + \sin\phi u_z^-).
\end{align}

Similarly, in the new coordinate system, the jump condition Eq.  (\ref{jumpcondition2}) can be discretized as
\begin{align}
\label{3dbasiccondition2}&[\beta u_{\xi}]=\beta^+ (\sin\phi \cos\theta u_x^+  + \sin \phi\sin \theta u_y^+ +\cos\phi u_z^+ ) - \beta^-(\sin\phi \cos \theta u_x^-  + \sin \phi\sin \theta u_y^- +\cos\phi u_z^- ).
\end{align}

These jump conditions, together with the original jump condition Eq.  (\ref{jumpcondition1}) , constitutes a set of lowest order jump conditions that is enforced in the MIB method.

In principle, by differentiating these lowest order jump conditions, we can generate even more jump conditions. However by doing this, we could create some higher-order derivatives and cross derivatives, whose evaluation often involves larger stencils, which is unstable for constructing high-order numerical schemes and unfeasible for truly complex interface geometries. Therefore, the spirit of the original MIB method is to use only the lowest order interface jump conditions. This spirit is preserved in  the proposed MIB-FVM. We  only use jump conditions in  Eqs.   (\ref{jumpcondition1}), (\ref{3dbasiccondition3}), (\ref{3dbasiccondition4}) and (\ref{3dbasiccondition2}).

Thus, at each intersection of the interface and the mesh line, there are four jump conditions described by Eqs.  (\ref{jumpcondition1}), (\ref{3dbasiccondition3}), (\ref{3dbasiccondition4}) and (\ref{3dbasiccondition2}), which involve  the function values and only the first order derivatives. This property is very important for making the numerical scheme more stable and efficient in handling complex interface geometries. Besides, employing only the lowest derivatives will generate a more-banded matrix with a smaller conditional number. The enforcement of these jump conditions determines fictitious values to be used in the discretization of differential operators in a given PDE.

The discretization of the differential operator involves six partial derivatives, namely $\frac{{\partial u}^{\pm}}{\partial x_j}$, where $x_j=x,y,z$. In dealing with complex interface geometries, some of these partial derivatives could be very difficult or even impossible to be approximated with some given order of accuracy. Since in constructing a second order numerical scheme, there is a pair of fictitious values along each mesh line and they are solved simultaneously, with the aforementioned four jump conditions, it is affordable for us to eliminate two of them. This redundancy gives two more degrees of freedom for us to design efficient and robust second order schemes in handling complex interface geometries. The designed numerical scheme  systematically eliminates the two partial derivatives that are most difficult to discretize at each intersection of the interface and the mesh line by using two redundant jump conditions. In summary, at each intersection of the interface and the mesh line, there are four remaining partial derivatives that are need to be discretized.

The elimination process must comply with the following rules:
\begin{itemize}
\item All fictitious values on all irregular points are determined. At the intersection of the interface and each mesh line, two fictitious values are determined at the same time.
\item Along each mesh line, two corresponding derivatives must be kept.
\item Among the remaining four partial derivatives, eliminate two that are most difficult to evaluate according to the local geometry.
\end{itemize}

From Eq.  (\ref{transformrelation}), we have
\begin{align}
	\left(\begin{array}{c}
	u_{\xi} \\
	u_{\eta} \\
	u_{\zeta}
	\end{array} \right) &={\bf P} \cdot \left(
	\begin{array}{c}
	u_x \\
	u_y \\
	u_z
	\end{array}\right)
\end{align}

Therefore, Eqs.   (\ref{3dbasiccondition2}), (\ref{3dbasiccondition3}) and (\ref{3dbasiccondition4}) can be rewritten as
\begin{align}\label{3ddiscretizedjumpcondition}
	\left(\begin{array}{c}
	{[\beta u_{\xi}]} \\
	{[u_{\eta}]} \\
	{[u_{\zeta}]}
	\end{array} \right)={\bf C} \cdot \left(
	\begin{array}{c}
	u_x^+ \\
	u_x^- \\
	u_y^+ \\
	u_y^- \\
	u_z^+ \\
	u_z^-
	\end{array} \right)
\end{align}
where
\begin{align}
&{\bf C}=\left(\begin{array}{c}
C_1 \\
C_2 \\
C_3
\end{array}\right)=\left(\begin{array}{cccccc}
p_{11}\beta^+ & -p_{11}\beta^- & p_{12}\beta^+ & -p_{12}\beta^- & p_{13}\beta^+ & -p_{13}\beta^- \\
p_{21} & -p_{21} & p_{22} & -p_{22} & p_{23} & -p_{23} \\
p_{31} & -p_{31} & p_{32} & -p_{32} & p_{33} & -p_{33}
\end{array}\right).
\end{align}
Here $p_{ij}$ is the $ij$th entry of the transformation matrix ${\bf C}$ and $C_i$ denotes $i$th row of ${\bf C}$. After eliminating the $l$th and the $m$th elements of entry of the array $(u_x^+,u_x^-,u_y^+,u_y^-,u_z^+,u_z^-)^T$, Eq.  (\ref{3ddiscretizedjumpcondition}) becomes
\begin{align} \label{aftereliminationjumpcondition}
a[\beta u_{\xi}]+b[u_{\eta}]+c[u_{\zeta}]=(aC_1+bC_2+cC_3) \cdot \left(\begin{array}{c}
	u_x^+ \\
	u_x^- \\
	u_y^+ \\
	u_y^- \\
	u_z^+ \\
	u_z^-
\end{array}\right)
\end{align}
where
\begin{align}
\label{eliminationcoefficients}
a&=C_{2l}C_{3m}-C_{3l}C_{2m} \\ \nonumber
b&=C_{3l}C_{lm}-C_{1l}C_{3m} \\ \nonumber
c&=C_{1l}C_{2m}-C_{2l}C_{lm}.
\end{align}
In practice, the two remaining jump conditions Eqs.   (\ref{jumpcondition1}) and (\ref{aftereliminationjumpcondition}) are used to determine a pair of fictitious values at the intersection point of the interface and each of the mesh line at a time.

\subsection{General MIB scheme}

\begin{figure}[!ht]
\begin{center}
\includegraphics[scale=0.3]{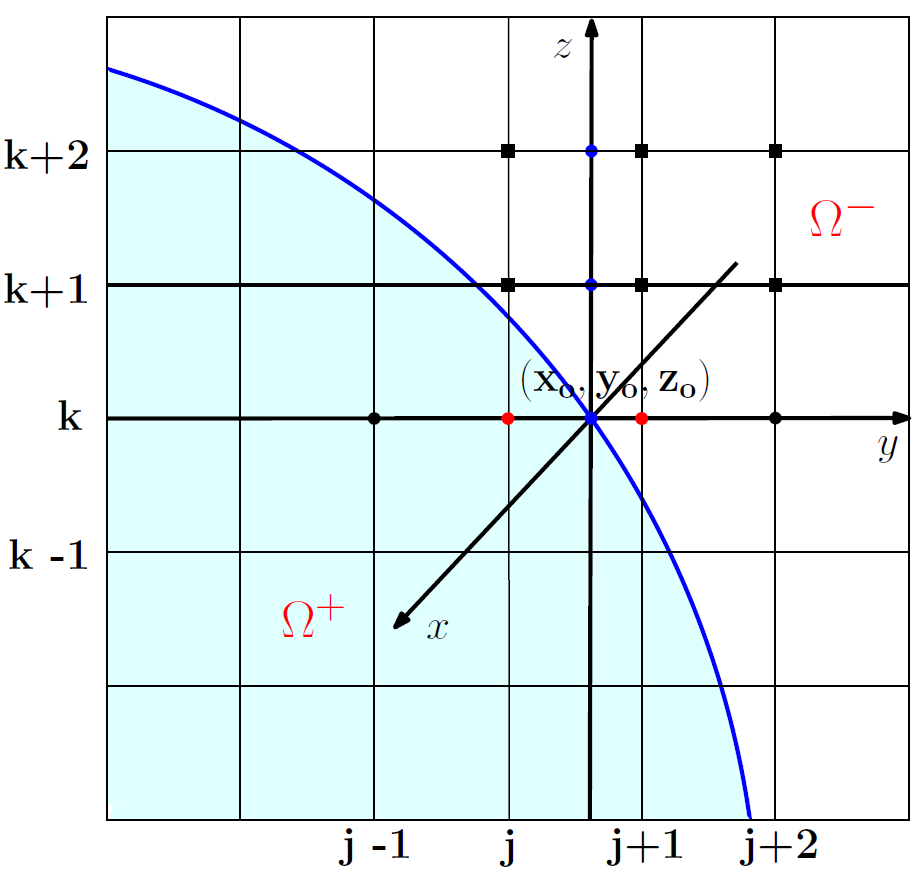}
\end{center}
\caption{Smooth interface at cross section $x=x_i$. The interface intersects with the $k$th mesh line in the $y$ direction at point $(x_o,y_o,z_o)$, where jump conditions need to be matched. At a pair of grid points in red color, namely $(i,j,k)$ and $(i,j+1,k)$, fictitious values are to be determined. Four grid points on the $y$ mesh line are employed to approximate quantities at $(x_o,y_o,z_o)$, and six grid points denoted by black squares are employed to approximate $u_z^-\mid_{(x_o.y_o,z_o)}$.}
\label{3dsmooth}
\end{figure}

Consider a geometry illustrated in Figure \ref{3dsmooth}. The domains are denoted by $\Omega^+$ and $\Omega^-$. The interface intersects the $k$th $y$-mesh line at point $(x_o,y_o,z_o)$ on the $yz$ plane. For discretizing central derivatives, two fictitious values, $f_{i,j,k}$ and $f_{i,j+1,k}$, are to be determined on irregular grid points.

Here, $u^+$, $u^-$, $u_y^+$ and $u_y^-$ at the intersection point $(x_o,y_0,z_0)$ are easily expressed by interpolation and the CFD scheme from information in $\Omega^+$ and $\Omega^-$, respectively
\begin{align}
u^+&=(w_{0,j-1},w_{0,j},w_{0,j+1}) \cdot (u_{i,j-1,k},u_{i,j,k},f_{i,j+1,k})^T \\
u^-&=(w_{0,j},w_{0,j+1},w_{0,j+2}) \cdot (f_{i,j,k},u_{i,j+1,k},u_{i,j+2,k})^T \\
u_y^+&=(w_{1,j-1},w_{1,j},w_{1,j+1}) \cdot (u_{i,j-1,k},u_{i,j,k},f_{i,j+1,k})^T \\
u_y^-&=(w_{1,j},w_{1,j+1},w_{1,j+2}) \cdot (f_{i,j,k},u_{i,j+1,k},u_{i,j+2,k})^T,
\end{align}
where $w_{s,t}$ denote interpolation weights, which is generated by using the standard Lagrange polynomials \cite{Fornberg:1998}. The first subscript $s$ represents either interpolation (when $s=0$) or first order derivative (when $s=1$) at point $(x_o,y_o,z_o)$, while the second subscript represents the node index.

We only need to compute two out of four remaining first order partial derivatives. For instance, if $u_x^-$ and $u_z^-$ are easily computed, then by choosing $s=1$ and $t=5$ in Eqs.   (\ref{aftereliminationjumpcondition}) and (\ref{eliminationcoefficients}), we can eliminate $u_x^+$ and $u_z^+$.

To approximate $u_z^+$ or $u_z^-$, we need three more function values at the intersection points of the auxiliary line $y=y_o$ and the mesh line on the $y-z$ plane. Here we only provide a detailed scheme to   approximate $u_z^-$, while other derivatives can be approximated in the same manner. Since these points are not on  grid nodes, they are interpolated by nearby grid points along the $y$ direction. Therefore, six more auxiliary points are involved. As illustrated in Figure \ref{3dsmooth}, $u_z^-(x_o,y_o,z_o)$ can be approximately as
{\small
\begin{align}
u_z^-(x_o,y_o,z_o)=(w_{1,k},w_{1,k+1},w_{1,k+2}) \cdot
\left(\begin{array}{ccccccccc}
w_{0,j} & w_{0,j+1} & w_{0,j+2} & 0 & 0 & 0 & 0 & 0 & 0 \\
0 & 0 & 0 & w_{0,j}' & w_{0,j+1}' & w_{0,j+2}' & 0 & 0 & 0 \\
0 & 0 & 0 & 0 & 0 & 0 & w_{0,j}^* & w_{0,j+1}^* & w_{0,j+2}^*
\end{array} \right) \cdot {\bf U}
\end{align}
}
where ${\bf U}=(f_{i,j,k},u_{i,j+1,k},u_{i,j+2,k},u_{i,j,k+1},u_{i,j+1,k+1},u_{i,j+2,k+1},u_{i,j,k+2},u_{i,j+1,k+2},u_{i,j+2,k+2})^T$. Here the superscripts on $w$ denotes different sets of FD weights \cite{Fornberg:1998}. Similarly, $u_x^-$ can be approximated along the auxiliary line $y=y_o$ on the $xy$ plane. Then jump conditions Eqs.   (\ref{jumpcondition1}) and (\ref{aftereliminationjumpcondition}) combined are used to solve the pair of fictitious values $f_{i,j,k}$ and $f_{i,j+1,k}$ at the same time. Consequently, the two fictitious values are expressed as a linear combination of 16 function values at nearby grid points and 4 jump conditions.

%

\subsection{MIB scheme for interface with large curvatures}

The above argument is based on a crucial assumption, which is that on each mesh line, there should be enough grid nodes (at least two for employing second order standard FD scheme) near the interface inside each subdomain so that the jump conditions can be expressed. However, when the curvature of the interface is very large, this assumption cannot be guaranteed for all mesh sizes.

\begin{figure}[!ht]
\begin{center}
\includegraphics[scale=0.3]{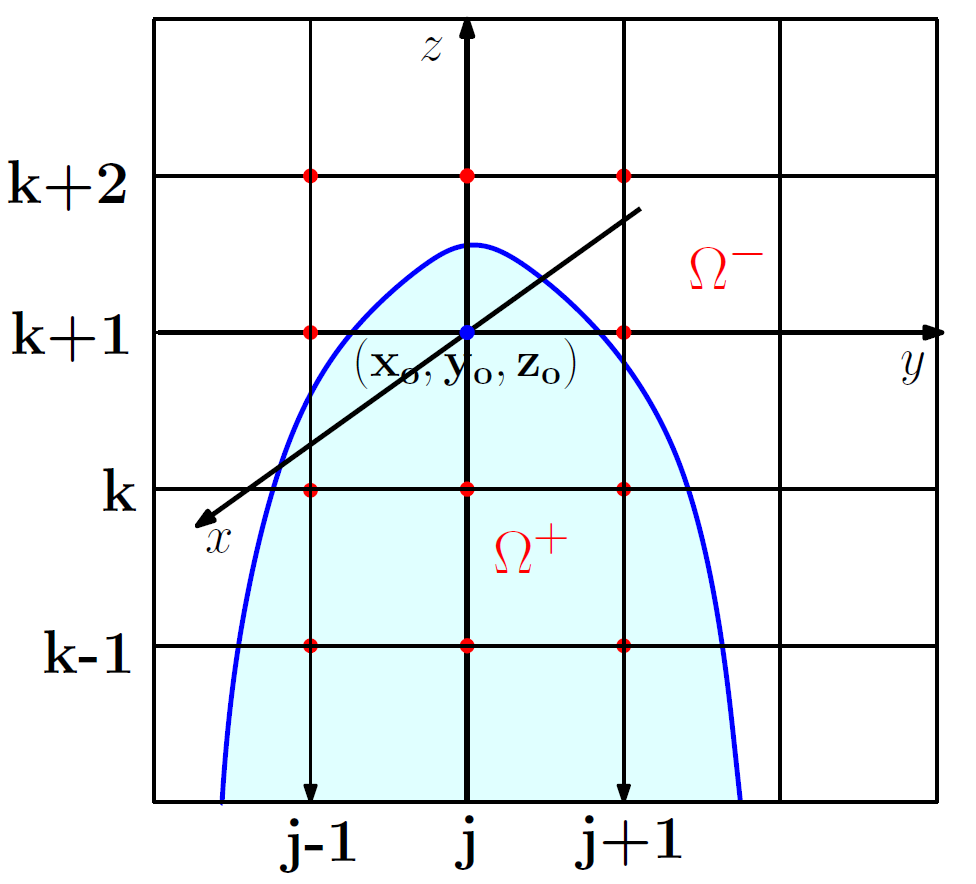}
\end{center}
\caption{Disassociation type of irregular grid points at cross section $x=x_i$. Fictitious value need to be deployed at point $(i,j,k)$, $(i,j-1,k)$ and $(i,j+1,k)$ in order to discretize the PDE at these points. However, $f_{i,j,k}$, $f_{i,j-1,k}$ and $f_{i,j+1,k}$ can not be obtained from the $y$-direction by the general MIB scheme, since there are not enough points in $\Omega^+$. Luckily, these fictitious value can be obtained from the $z$-direction. For the above figure, during the discretization process, fictitious values obtained in the $z$-direction will be utilized for both the $y$-direction and the $z$-direction discretizations of the PDE.}
\label{dissassociation}
\end{figure}

As is illustrated in Figure \ref{dissassociation}, fictitious value $f_{i,j,k}$, $f_{i,j-1,k}$ and $f_{i,j+1,k}$ cannot be solved along $y$-direction, since there is only one grid point inside the interface in the $y$-direction. The concept of disassociation is introduced in Ref. \cite{Zhou:2006d} in order to disassociate the domain extension from the discretization thus broaden the applicability of the MIB method to general interface geometry. In the situation of Figure \ref{dissassociation}, we cannot solve the fictitious values in the $y$-direction, however, there is no problem to solve the fictitious values in the $z$-direction. Consequently, the fictitious value at an irregular point, regardless of the direction in its calculation, can be used for discretization in any direction involving the grid point without loss of accuracy. For fictitious values that can be obtained in multiple directions, we may compute their values in the most convenient manner and use them for necessary discretization. In practice, in order to make the MIB matrix as symmetric and banded as possible, if the fictitious value can be found in the given direction, one should use the fictitious value obtained from the corresponding direction and avoid using the disassociation technique. It is also worth of mentioning that disassociation technique does not reduce the accuracy of the numerical scheme. For instance, if fictitious values obtained from the $z$-direction has $O(h^m)$ accuracy for some integer $m$. This accuracy of approximation is independent of the direction in which the fictitious values are obtained. Here $h$ is the grid size of the uniform mesh.

\newpage
\section{Numerical studies}\label{Numericalstudies}

In this section, we examine validity and test the performance of the proposed MIB-FVM for solving the 3D Poisson equation with discontinuous coefficients. To demonstrate the robustness and accuracy of the present MIB-FVM, we carry out many numerical experiments with complex interfaces in both 2D and 3D settings. For 2D test cases, we start with a 6-petal flower shape interface which is relatively complex. The 2D jigsaw puzzle-like interface is also studied. For 3D  problems, we consider various interface geometries including  sphere,  ellipsoid, cylinder, 5-petal flower and torus. Finally, we investigate our MIB-FVM's ability to deal  with low regularity solutions in two more test cases.

{\bf Case 1.}
We first consider a 2D interface problem. The 2D Poisson equation is solved in domain $[-1,1]\times[-1,1]$. The interface of the problem is a 6-petal flower whose formula in polar coordinate system is given as $r=\frac{1}{2}(1+\frac{1}{2}\sin(6\theta))$.
The designed solution in two different domains are given by
\begin{align}
u^+ (x,y)=\frac{1}{4}+\sin(x)\sin(y),  \quad {\rm and} \quad  u^- (x,y)=0.
\end{align}
The coefficients are given by $\beta^+ (x,y)=\beta^- (x,y)=1$.

\begin{figure}
\begin{center}
\begin{tabular}{cc}
\includegraphics[width=0.5\textwidth]{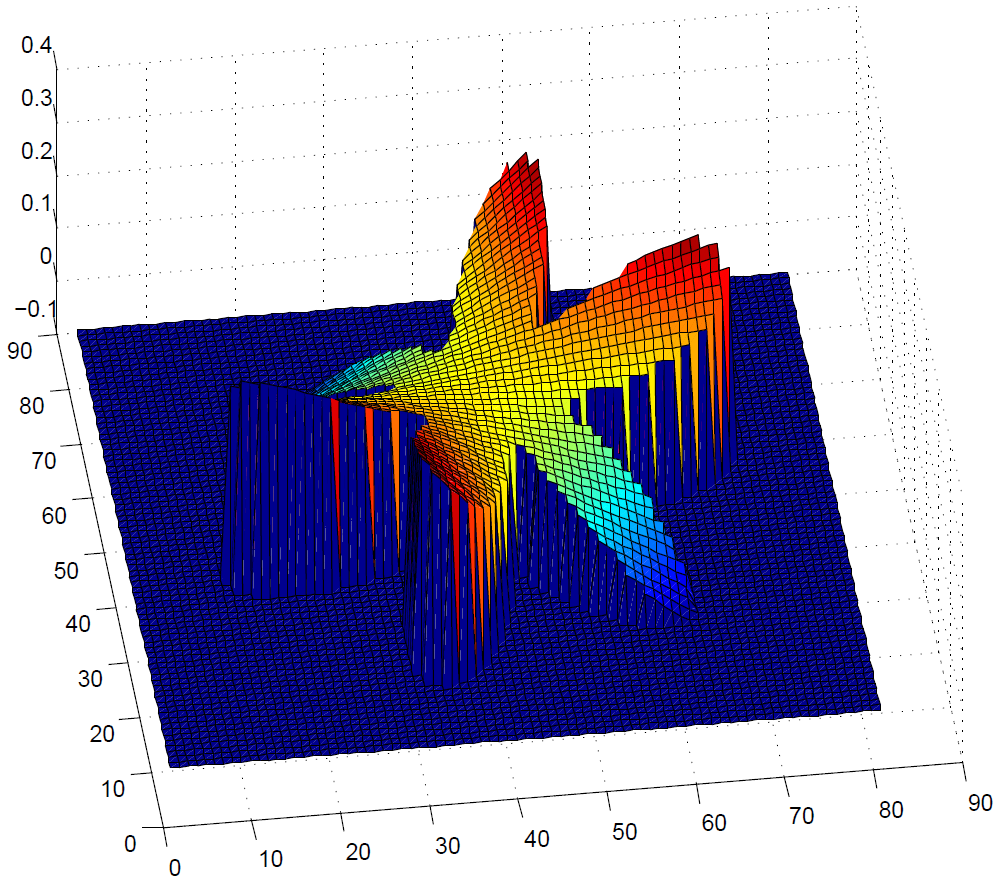} &
\includegraphics[width=0.5\textwidth]{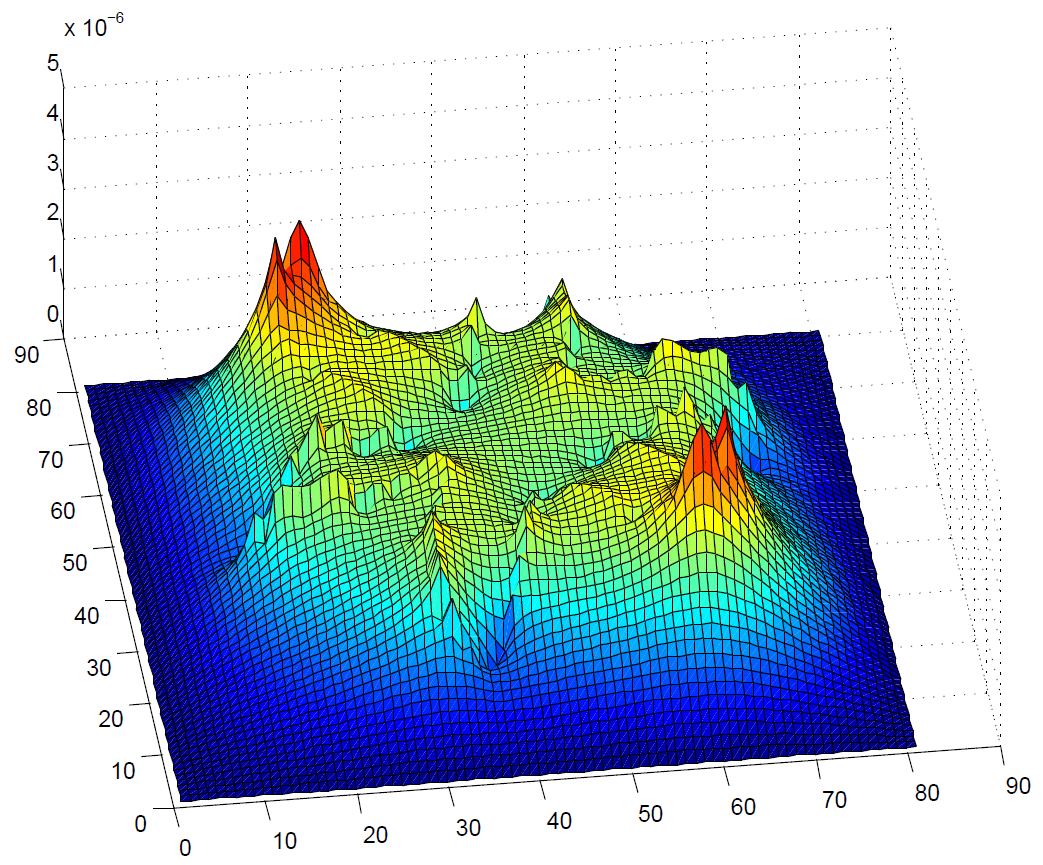}
\end{tabular}
\end{center}
\caption{The computed solution on an $80 \times 80$  mesh (Left chart) and the $L_{\infty}$  error (Right chart) for Case 1.}
\label{6petal_pic}
\end{figure}

\begin{table}
\caption{Numerical errors and convergence orders for a 6-petal flower interface problem (Case 1).}
\label{6petal_tab}
\begin{center}
\begin{tabular}{c||cc|cc}
\hline
$n_x \times n_y$ &$L_{\infty}(u)$ &Order   &$L_2 (u )$ &Order     \\\hline
40 $\times$ 40   & 3.790e-5       &       &1.910e-5   &            \\
80 $\times$ 80   & 4.867e-6       & 2.96  &1.827e-6   & 3.38        \\
160 $\times$ 160 & 8.545e-7       & 2.51  &2.275e-7   & 2.99        \\
\hline
\end{tabular}
\end{center}
\end{table}

Figure \ref{6petal_pic} depicted the numerical solution and the $L_{\infty}$ error on a $40 \times 40 \times 40$ mesh. From the plot of the $L_{\infty}$ error on the right, it is noticed that the highest error appears at the tips of the petals, where the curvature of the interface is relatively higher. Table \ref{6petal_tab} shows the results of the numerical accuracy tests on three successively refined meshes. From the table, it can be seen that  MIB-FVM    obtains the second order accuracy in both the $L_{\infty}$ error and $L_2$ error.


{\bf Case 2.}
We next consider the 2D jigsaw puzzle-like interface problem. The 2D Poisson equation is solved in domain $[-1,1]\times[0,3]$. The interface of the problem is a jigsaw puzzle-like shape whose formula is given as
\begin{align}
	\begin{cases}
	x(\theta)=0.6\cos(\theta)-0.3\cos(3\theta) \\
	y(\theta)=1.5+0.7\sin(\theta)-0.07\sin(3\theta)+0.2\sin(7\theta)
	\end{cases}
\end{align}
The solution in two different subdomains is chosen as
\begin{align}
u^+ (x,y)=e^x(y^2+x^2\sin(y)),  \quad  {\rm and} \quad  u^- (x,y)=-(x^2+y^2)
\end{align}
The discontinuous coefficients are given by $\beta^+ (x,y,z)=1, \quad {\rm and} \quad \beta^- (x,y,z)=10$.

\begin{figure}
\begin{center}
\begin{tabular}{ccc}
\includegraphics[width=0.5\textwidth]{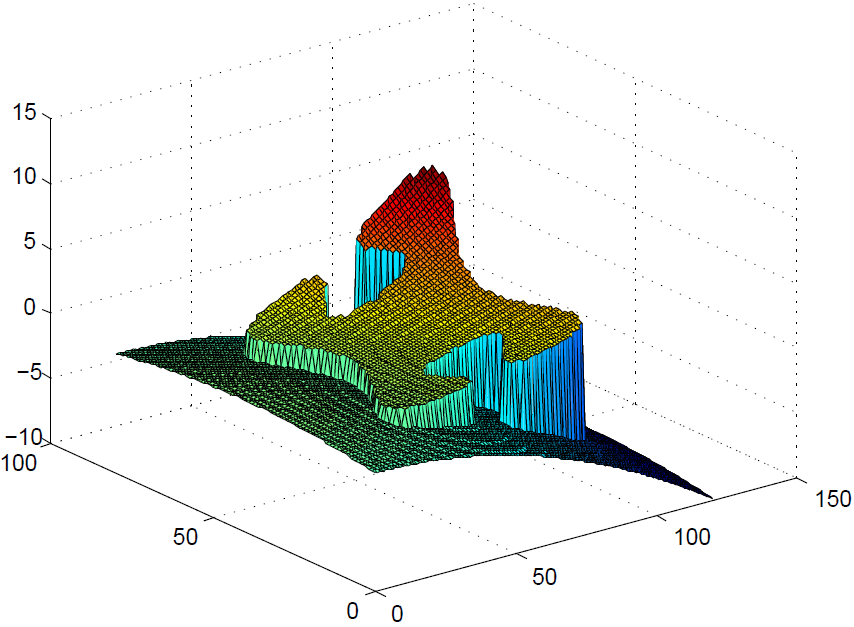} &&
\includegraphics[width=0.5\textwidth]{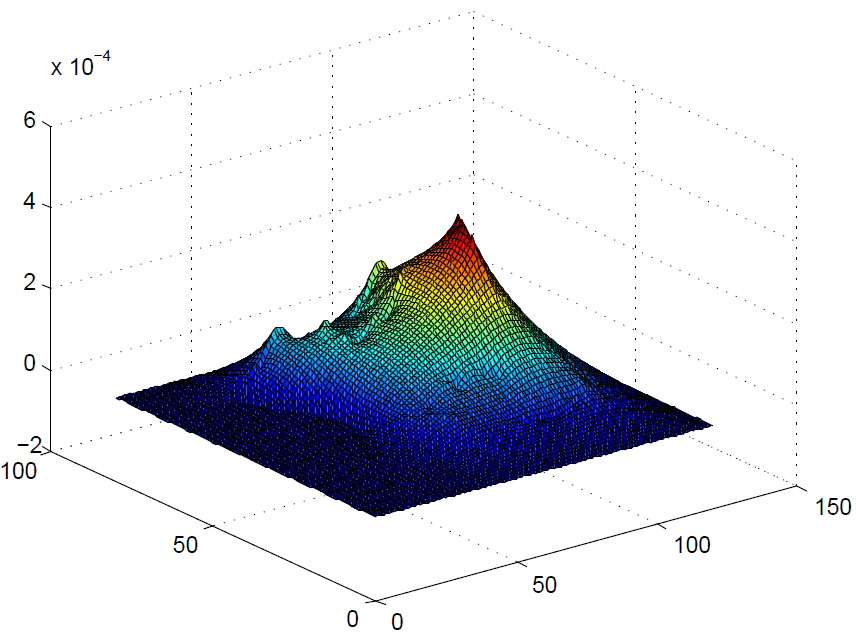}
\end{tabular}
\end{center}
\caption{The computed solution on a  $40 \times 40$  mesh (Left chart) and the $L_{\infty}$  error (Right chart) for Case 2.}
\label{puzzle_pic}
\end{figure}

\begin{table}
\caption{Numerical  errors and convergence orders for an ellipsoidal interface problem (Case 2).}
\label{puzzle_tab}
\begin{center}
\begin{tabular}{c||cc|cc}
\hline
$n_x \times n_y$  &$L_{\infty}(u)$ &Order   &$L_2 (u )$ &Order     \\\hline
40 $\times$ 60    & 1.490e-3  &       &2.176e-4   &            \\
80 $\times$ 120   & 2.255e-4  & 2.72  &3.654e-5   & 2.57        \\
160$\times$ 240   & 2.448e-5  & 3.20  &4.067e-6   & 3.16         \\
\hline
\end{tabular}
\end{center}
\end{table}

Figure \ref{puzzle_pic} depicted the computed solution and the $L_{\infty}$ error on a $40 \times 40 \times 40$ mesh. The numerical errors in terms of $L_{\infty}$ norm are collected in Table \ref{puzzle_tab}, which also show  the second order convergence of the proposed MIB-FVM.


{\bf Case 3.}
Next we investigate a classical spherical interface problem. The 3D Poisson equation is solved in domain $[-1,1]\times[-1,1]\times[-1,1]$. To specify the spherical interface, we design the level set function $\phi(x,y,z)$
\begin{align}
\phi(x,y,z)= r_0^2-(x^2+y^2+z^2),
\end{align}
where the radius $r_0=3/4$. The interface is given as $\Gamma=\{(x,y,z)| \phi=0, \forall (x,y,z)\in \Omega \}$. Two subdomains are $\Omega^+=\{(x,y,z)| \phi \geq0, \forall  (x,y,z)\in \Omega \}$, and $\Omega^-=\{(x,y,z)| \phi < 0; \forall  (x,y,z)\in \Omega\}$.
The solution in two different subdomains is chosen as
\begin{align}
u^+ (x,y,z)=0,  \quad {\rm and} \quad  u^- (x,y,z)=\cos(x)\cos(y)\cos(z).
\end{align}
The discontinuous coefficients are given by
\begin{itemize}
\item {\bf Case 3(a):}
\begin{align}
\beta^+ (x,y,z)=8, \quad {\rm and} \quad \beta^- (x,y,z)=1.
\end{align}
\item {\bf Case 3(b):}
\begin{align}
\beta^+ (x,y,z)=1+\cos\left(\frac{x}{10}+\frac{y}{10}+\frac{z}{10}\right), \quad {\rm and} \quad \beta^- (x,y,z)=1+3\sin\left(\frac{x}{10}+\frac{y}{10}+\frac{z}{10}\right).
\end{align}
\end{itemize}

\begin{figure}
\begin{center}
\begin{tabular}{ccc}
\includegraphics[width=0.49\textwidth]{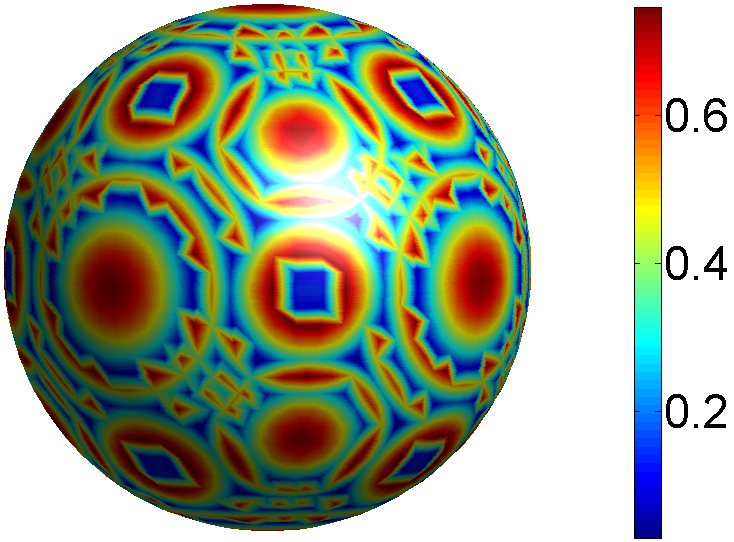} &&
\includegraphics[width=0.51\textwidth]{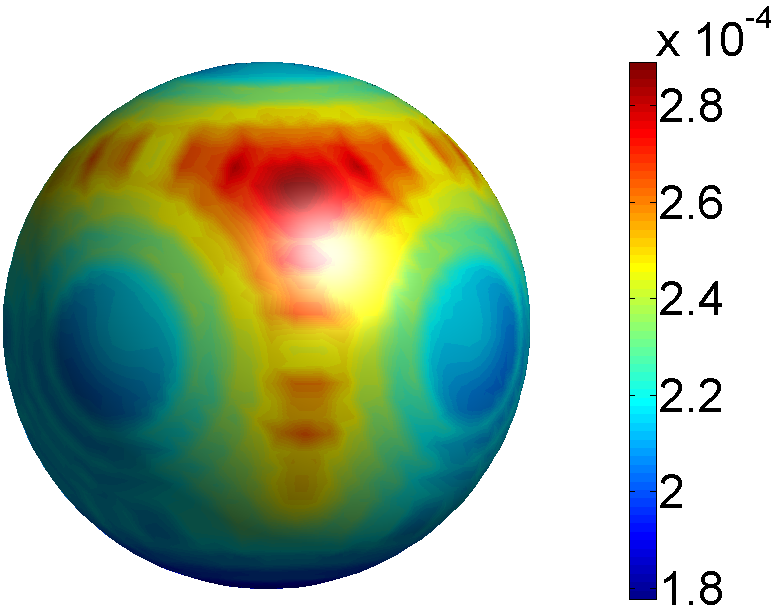}
\end{tabular}
\end{center}
\caption{The computed solution on a  $40 \times 40 \times 40$  mesh (Left chart) and the $L_{\infty}$  error (Right chart) for Case 3(a).}
\label{sphere1_pic}
\end{figure}

\begin{table}
\caption{Numerical errors and convergence orders for a spherical interface problem (Case 3(a)).}
\label{sphere1_tab}
\begin{center}
\begin{tabular}{c||cc|cc}
\hline
$n_x \times n_y \times n_z$ &$L_{\infty}(u)$ &Order   &$L_2 (u )$ &Order     \\\hline
20 $\times$ 20 $\times$ 20   & 1.2771e-3  &       &5.5733e-4   &            \\
40 $\times$ 40 $\times$ 40   & 2.8944e-4  & 2.14  &1.2612e-4   & 2.14        \\
80 $\times$ 80 $\times$ 80   & 6.4798e-6  & 2.15  &2.8204e-5   & 2.16        \\
160 $\times$ 160 $\times$ 160 & 1.1673e-5  & 2.47  &4.6283e-6   & 2.61
\\\hline
\end{tabular}
\end{center}
\end{table}

Table \ref{sphere1_tab} lists numerical results for the spherical interface problem. The geometry and error are depicted in Figure \ref{sphere1_pic}. It is seen that errors are very small at the fine mesh. The second order convergence is found.

To further investigate the present method, we change the coefficient into a position dependent function in Case 3(b). Result is given in Table \ref{sphere2_tab}, which also demonstrates the second order accuracy in both $L_{\infty}$ and $L_2$ norms for the solution.

\begin{figure}
\begin{center}
\begin{tabular}{ccc}
\includegraphics[width=0.49\textwidth]{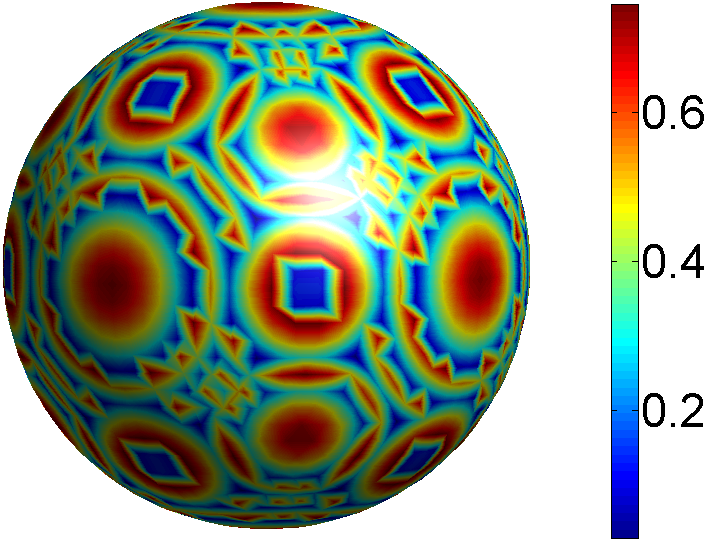} &&
\includegraphics[width=0.51\textwidth]{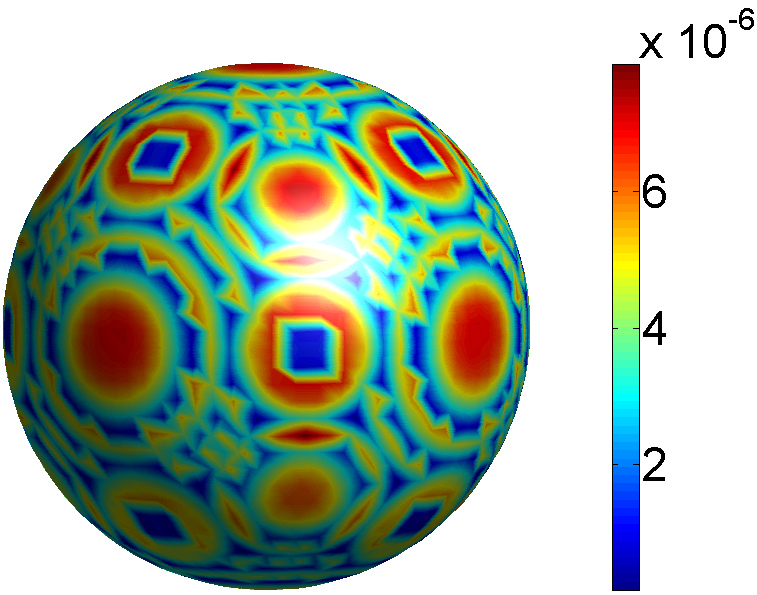}
\end{tabular}
\end{center}
\caption{The computed solution on a  $40 \times 40 \times 40$  mesh (Left chart) and the $L_{\infty}$  error (Right chart) for Case 3(b).}
\label{sphere2_pic}
\end{figure}

\begin{table}
\caption{Numerical  errors and convergence orders for a spherical interface problem (Case 3(b)).}
\label{sphere2_tab}
\begin{center}
\begin{tabular}{c||cc|cc}
\hline
$n_x \times n_y \times n_z$ &$L_{\infty}(u)$ &Order   &$L_2 (u )$ &Order     \\\hline
20 $\times$ 20 $\times$ 20   & 8.3466e-5  &       &3.0069e-5   &            \\
40 $\times$ 40 $\times$ 40   & 1.9207e-5  & 2.12  &6.8250e-6   & 2.14        \\
80 $\times$ 80 $\times$ 80   & 4.9351e-6  & 1.96  &1.5630e-6   & 2.13        \\
160 $\times$ 160 $\times$ 160   & 1.5551e-6  & 1.67  &3.7873e-7   & 2.05        \\
\hline
\end{tabular}
\end{center}
\end{table}


{\bf Case 4.}
We next consider an  ellipsoidal interface problem. The 3D Poisson equation is solved in domain $[-5,5]\times[-5,5]\times[-5,5]$. The interface of the problem is an ellipsoid given as $\left(\frac{x}{\frac{2}{7}}\right)^2+\left(\frac{y}{\frac{25}{14}}\right)^2+\left(\frac{z}{\frac{25}{14}}\right)^2=1$. To specify the ellipsoidal interface, we define the level set function $\phi(x,y,z)$
\begin{align}
\phi(x,y,z)=1-\left(\left(\frac{x}{\frac{2}{7}}\right)^2+\left(\frac{y}{\frac{25}{14}}\right)^2+\left(\frac{z}{\frac{25}{14}}\right)^2\right)
\end{align}
The interface is given as $\Gamma=\{(x,y,z)| \phi=0, \forall (x,y,z)\in \Omega \}$. Two subdomains are $\Omega^+=\{(x,y,z)| \phi \geq0, \forall  (x,y,z)\in \Omega \}$, and $\Omega^-=\{(x,y,z)| \phi < 0; \forall  (x,y,z)\in \Omega\}$.
The solution in two different subdomains is chosen as
\begin{align}
u^+ (x,y,z)=e^{-\frac{x^2+y^2+z^2}{20}},  \quad {\rm and} \quad  u^- (x,y,z)=\cos(x)\cos(y)\cos(z).
\end{align}
The discontinuous coefficients are given by $\beta^+ (x,y,z)=z+15, \quad {\rm and} \quad \beta^- (x,y,z)=\frac{x+y}{2}+10$.

\begin{figure}
\begin{center}
\begin{tabular}{ccc}
\includegraphics[width=0.495\textwidth]{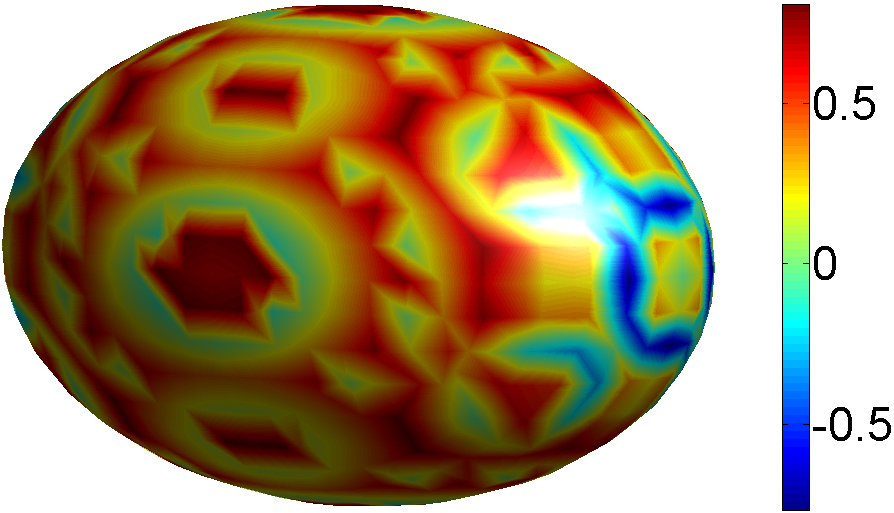} &&
\includegraphics[width=0.505\textwidth]{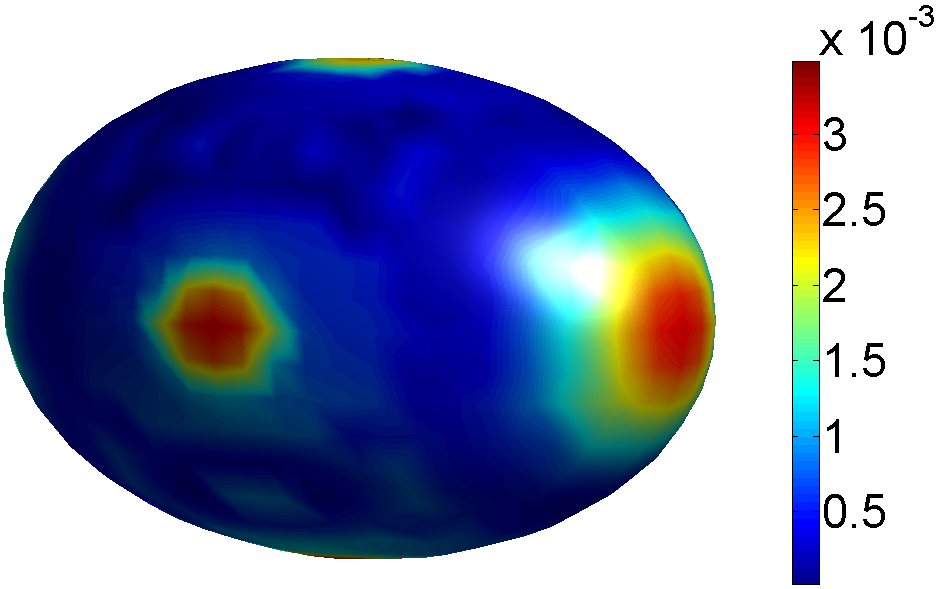}
\end{tabular}
\end{center}
\caption{The computed solution on a  $40 \times 40 \times 40$  mesh (Left chart) and the $L_{\infty}$  error (Right chart) for Case 4.}
\label{ellipsoid_pic}
\end{figure}

\begin{table}
\caption{Numerical  errors and convergence orders for an ellipsoidal interface problem (Case 4).}
\label{ellipsoid_tab}
\begin{center}
\begin{tabular}{c||cc|cc}
\hline
$n_x \times n_y \times n_z$ &$L_{\infty}(u)$ &Order   &$L_2 (u )$ &Order     \\\hline
20 $\times$ 20 $\times$ 20   & 2.1881e-2  &       &6.8115e-3   &            \\
40 $\times$ 40 $\times$ 40   & 5.5719e-3  & 1.97  &1.7463e-3   & 1.96        \\
80 $\times$ 80 $\times$ 80   & 1.4675e-3  & 1.92  &4.5309e-4   & 1.95        \\
160 $\times$ 160 $\times$ 160   & 3.6068e-4  & 2.02  &1.0974e-4   & 2.05        \\
\hline
\end{tabular}
\end{center}
\end{table}

Figure \ref{ellipsoid_pic} depicted the numerical solution and the $L_{\infty}$ error on a $40 \times 40 \times 40$ mesh. From the result demonstrated in Table \ref{ellipsoid_tab}, it can be seen that the MIB-FVM   obtains the second order accuracy.


{\bf Case 5.}
In Case 5, the 3D Poisson equation is solved in domain $[-4,4]\times[-4,4]\times[-2,8.4]$. The interface of the problem is a cylinder of height $2\pi$ and base radius $\pi$. The designed solution in two different subdomains is chosen as
\begin{align}
u^+ (x,y,z)=x+y+z,  \quad {\rm and} \quad  u^- (x,y,z)=\cos(x)\cos(y)\cos(z).
\end{align}
The discontinuous coefficients are given by $\beta^+ (x,y,z)=8, \quad {\rm and} \quad \beta^- (x,y,z)=1$.

\begin{figure}
\begin{center}
\begin{tabular}{ccc}
\includegraphics[width=0.475\textwidth]{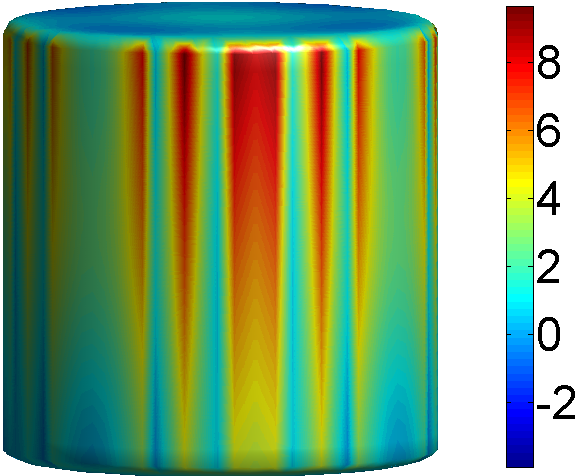} &&
\includegraphics[width=0.525\textwidth]{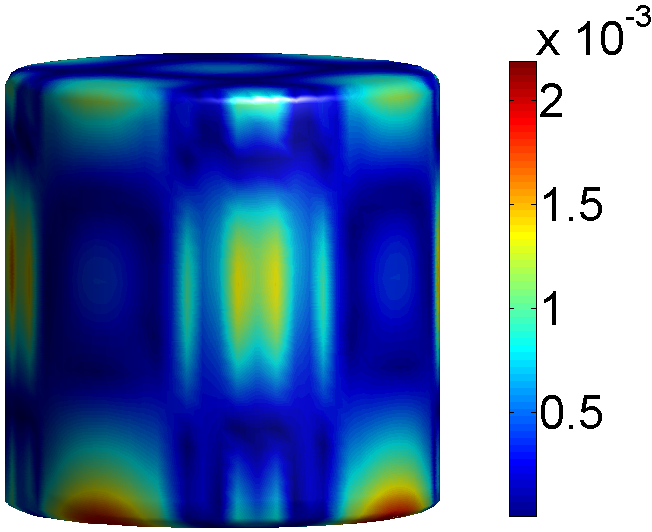}
\end{tabular}
\end{center}
\caption{The computed solution on a $40 \times 40 \times 40$  mesh (Left chart) and the $L_{\infty}$  error (Right chart) for Case 5.}
\label{cylinder_pic}
\end{figure}

\begin{table}
\caption{Numerical  errors and convergence orders for a cylinder interface problem (Case 5).}
\label{cylinder_tab}
\begin{center}
\begin{tabular}{c||cc|cc}
\hline
$n_x \times n_y \times n_z$ &$L_{\infty}(u)$ &Order   &$L_2 (u )$ &Order     \\\hline
20 $\times$ 20 $\times$ 20   & 1.2274e-2  &       &2.9195e-3   &            \\
40 $\times$ 40 $\times$ 40   & 3.0505e-3  & 2.01  &7.7935e-4   & 1.91        \\
80 $\times$ 80 $\times$ 80   & 7.9267e-4  & 1.94  &2.0674e-4   & 1.91        \\
160 $\times$ 160 $\times$ 160   & 1.8598e-4  & 2.09  &4.6819e-5   & 2.14      \\
\hline
\end{tabular}
\end{center}
\end{table}

Figure \ref{cylinder_pic} depicted the numerical solution and the $L_{\infty}$ error on a $40 \times 40 \times 40$ mesh. From the result demonstrated in Table \ref{cylinder_tab}, it can be seen that the proposed  MIB-FVM     achieves the second order accuracy.


{\bf Case 6.}
In Case 6, the 3D Poisson equation is solved in domain $[-5,5]\times[-5,5]\times[-2,2]$. The interface of the problem is a 5-petal flower-base cylinder whose polar equation is given by $r=\frac{5}{2}+\frac{5}{7}\sin(5\theta)$ and $-\frac{2}{3}\leq z \leq \frac{2}{3}$.
 The designed solution in two different subdomains is chosen as
\begin{align}
u^+ (x,y,z)=\sin(x)\sin(y)\sin(z),  \quad {\rm and} \quad  u^- (x,y,z)=e^{\frac{x}{10}+\frac{y}{10}+\frac{z}{10}}.
\end{align}
The discontinuous coefficients are given by $\beta^+ (x,y,z)=8, \quad {\rm and} \quad \beta^- (x,y,z)=1$.

\begin{figure}
\begin{center}
\begin{tabular}{ccc}
\includegraphics[width=0.49\textwidth]{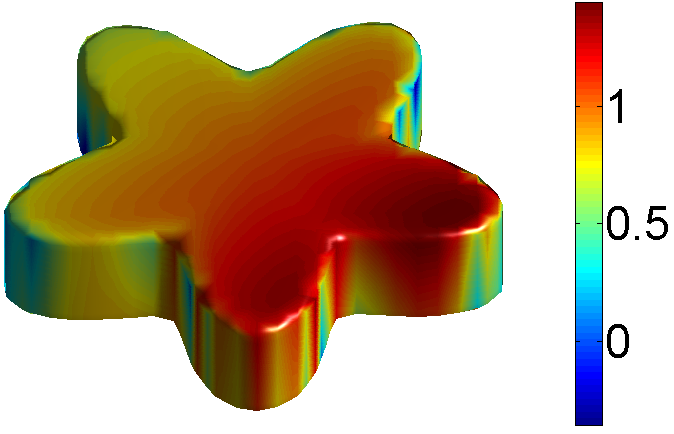} &&
\includegraphics[width=0.51\textwidth]{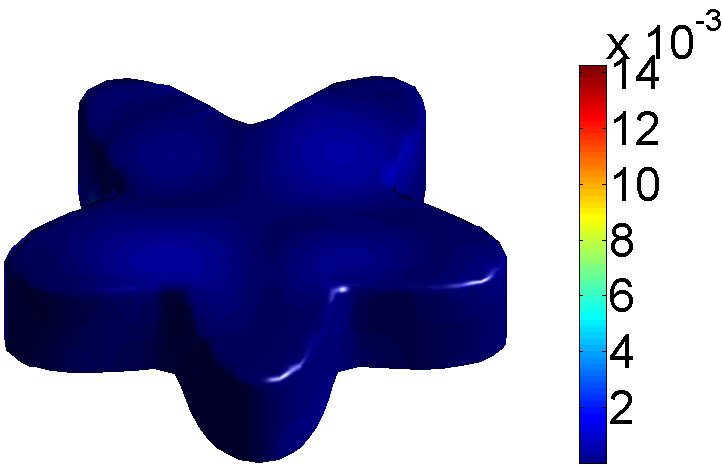}
\end{tabular}
\end{center}
\caption{The computed solution on a $40 \times 40 \times 40$  mesh (Left chart) and the $L_{\infty}$  error (Right chart) for Case 6.}
\label{fcylinder_pic}
\end{figure}

\begin{table}
\caption{Numerical errors and convergence orders for a 5-petal flower-base cylinder interface problem (Case 6).}
\label{fcylinder_tab}
\begin{center}
\begin{tabular}{c||cc|cc}
\hline
$n_x \times n_y \times n_z$ &$L_{\infty}(u)$ &Order   &$L_2 (u )$ &Order     \\\hline
20 $\times$ 20 $\times$ 20   & 4.9246e-2  &     &1.4819e-3   &         \\
40 $\times$ 40 $\times$ 40   & 1.4307e-2  & 1.78    &3.1350e-4   & 2.24        \\
80 $\times$ 80 $\times$ 80   & 4.9600e-4  & 4.85  &5.6857e-5   & 2.46        \\
160 $\times$ 160 $\times$ 160  &1.2096e-4  & 2.04  & 1.5637e-5   & 1.86           \\
\hline
\end{tabular}
\end{center}
\end{table}

Figure \ref{fcylinder_pic} depicted the numerical solution and the $L_{\infty}$ error on a $40 \times 40 \times 40$ mesh. From the result demonstrated in Table \ref{fcylinder_tab}, it can be seen that the present MIB-FVM  obtains the second order accuracy.


{\bf Case 7.}
In Case 7, the 3D Poisson equation is solved in domain $[-5,5]\times[-5,5]\times[-2,2]$. The interface of the problem is a torus whose equation is given by $(3-\sqrt{x^2+y^2})^2+z^2=1$.
The designed solution in two different subdomains is chosen as
\begin{align}
u^+ (x,y,z)=\sin(x)\sin(y)\sin(z)+1,  \quad {\rm and} \quad  u^- (x,y,z)=\cos(x)\cos(y)\cos(z).
\end{align}
The discontinuous coefficients are given by $\beta^+ (x,y,z)=8, \quad {\rm and} \quad \beta^- (x,y,z)=1$.

\begin{figure}
\begin{center}
\begin{tabular}{ccc}
\includegraphics[width=0.5\textwidth]{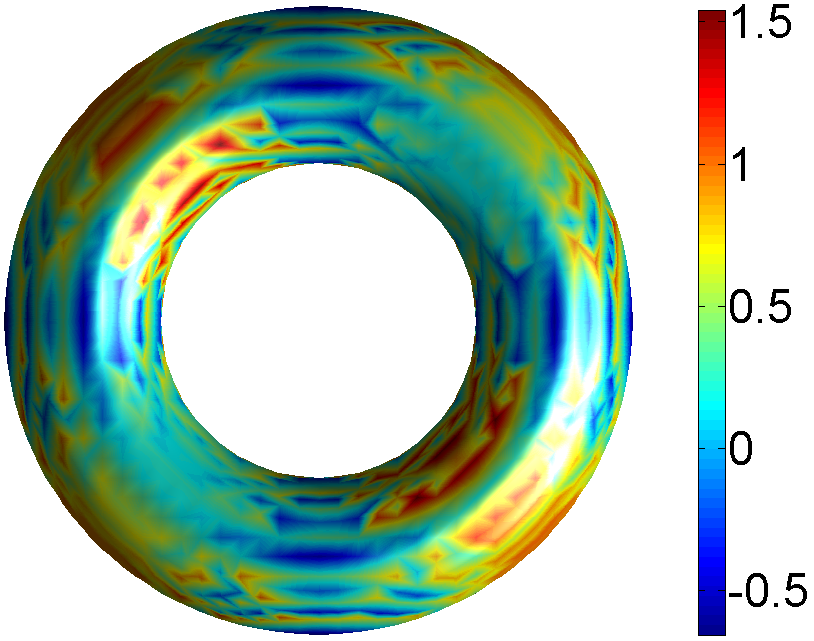} &&
\includegraphics[width=0.5\textwidth]{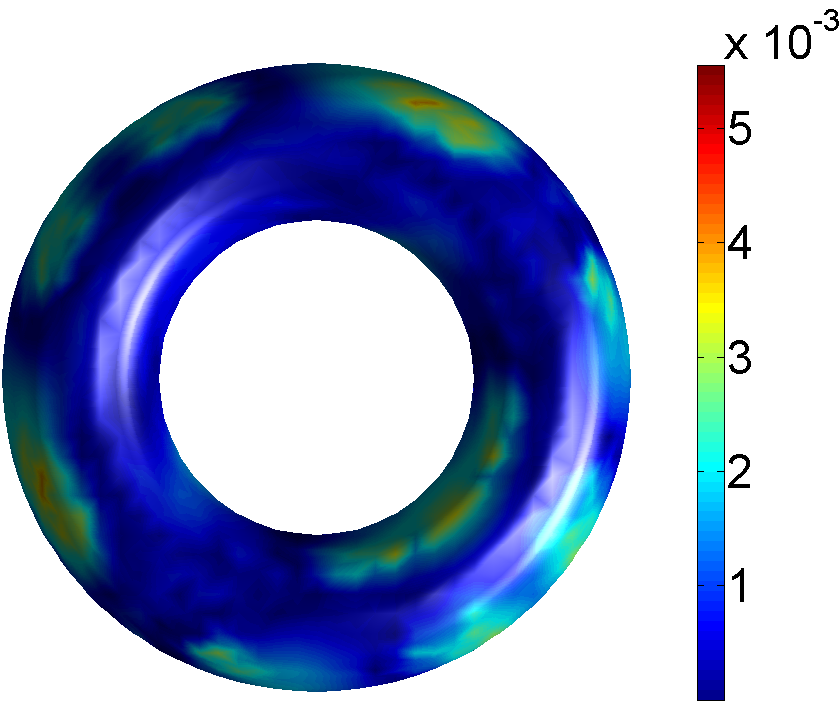}
\end{tabular}
\end{center}
\caption{The computed solution on a $40 \times 40 \times 40$  mesh (Left chart) and the $L_{\infty}$  error (Right chart) for Case 7.}
\label{torus_pic}
\end{figure}

\begin{table}
\caption{Numerical  errors and convergence orders for a torus interface problem (Case 7).}
\label{torus_tab}
\begin{center}
\begin{tabular}{c||cc|cc}
\hline
$n_x \times n_y \times n_z$ &$L_{\infty}(u)$ &Order   &$L_2 (u )$ &Order     \\\hline
20 $\times$ 20 $\times$ 20   & 3.4040e-2  &       &4.6519e-3   &            \\
40 $\times$ 40 $\times$ 40   & 5.7226e-3  & 2.57  &1.0688e-3   & 2.12        \\
80 $\times$ 80 $\times$ 80   & 1.6415e-3  & 1.80  &2.6468e-4   & 2.01         \\
160 $\times$ 160 $\times$ 160   & 3.8950e-4  & 2.08  &6.9023e-5   & 1.94         \\
\hline
\end{tabular}
\end{center}
\end{table}

Figure \ref{torus_pic} depicted the numerical solution and the $L_{\infty}$ error on a $40 \times 40 \times 40$ mesh. From the result demonstrated in Table \ref{torus_tab}, it can be seen that the present MIB-FVM    attains the second order accuracy.

{\bf Case 8.}
in this case, we test the proposed method for problems with low solution regularity.
\begin{itemize}
\item {\bf Case 8(a):}
The 3D Poisson equation is solved in domain $[-1,1.05]\times[-1,1.05]\times[-1,1.05]$. The interface of the problem is a sphere the same as that in Case 3.
The designed solution in two different subdomains is chosen as
\begin{align}
u^+ (x,y,z)=8,  \quad {\rm and} \quad  u^- (x,y,z)=(x^2+y^2+z^2)^{\frac{5}{6}}.
\end{align}
The discontinuous coefficients are given by $\beta^+ (x,y,z)=4, \quad {\rm and} \quad \beta^- (x,y,z)=1$.
\item {\bf Case 8(b):}
The 3D Poisson equation is solved in domain $[-5,5.05]\times[-5,5.05]\times[-5,5.05]$. The interface of the problem is an ellipsoid the same as that in Case 4.
The designed solution in two different subdomains is chosen as
\begin{align}
u^+ (x,y,z)=8,  \quad {\rm and} \quad  u^- (x,y,z)=(x^2+y^2+z^2)^{\frac{5}{6}}+\sin(x+y+z).
\end{align}
The discontinuous coefficients are given by $\beta^+ (x,y,z)=z+5, \quad {\rm and} \quad \beta^- (x,y,z)=\frac{x+y}{2}+10$.
\end{itemize}

\begin{figure}
\begin{center}
\begin{tabular}{ccc}
\includegraphics[width=0.48\textwidth]{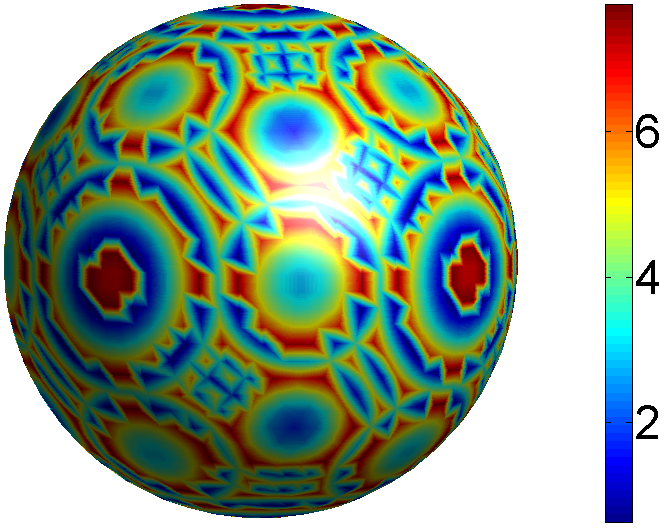} &&
\includegraphics[width=0.52\textwidth]{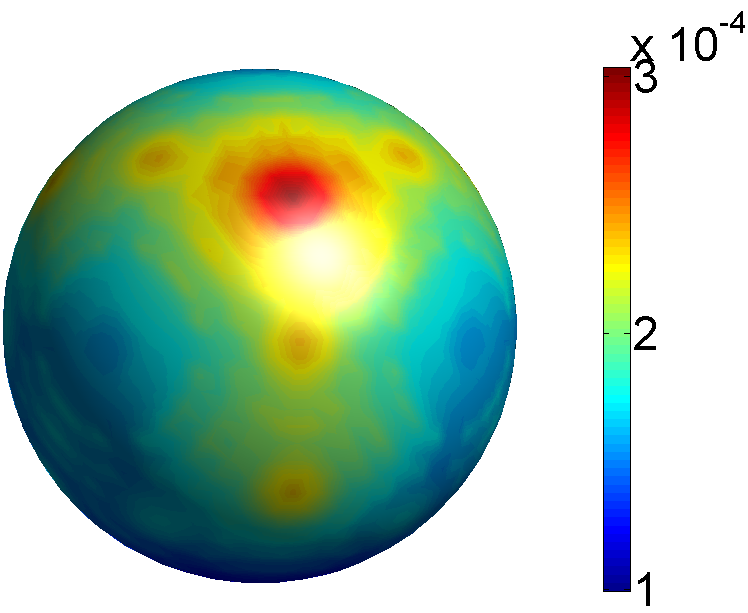}
\end{tabular}
\end{center}
\caption{The computed solution on a  $40 \times 40 \times 40$  mesh (Left chart) and the $L_{\infty}$  error (Right chart) for Case 8(a).}
\label{h2sphere_pic}
\end{figure}

\begin{table}
\caption{Numerical errors and convergence orders for a sphere interface problem (Case 8(a)).}
\label{h2sphere_tab}
\begin{center}
\begin{tabular}{c||cc|cc}
\hline
$n_x \times n_y \times n_z$ &$L_{\infty}(u)$ &Order   &$L_2 (u )$ &Order     \\\hline
20 $\times$ 20 $\times$ 20   & 1.1578e-3  &       &3.3263e-4   &            \\
40 $\times$ 40 $\times$ 40   & 3.0664e-4  & 1.92  &9.2064e-5   & 1.85        \\
80 $\times$ 80 $\times$ 80   & 6.8400e-5  & 2.17  &2.0947e-5   & 2.14         \\
160 $\times$ 160 $\times$ 160   & 1.3566e-5  & 2.34  &4.4747e-6   & 2.24       \\
\hline
\end{tabular}
\end{center}
\end{table}


\begin{figure}
\begin{center}
\begin{tabular}{cc}
\includegraphics[width=0.49\textwidth]{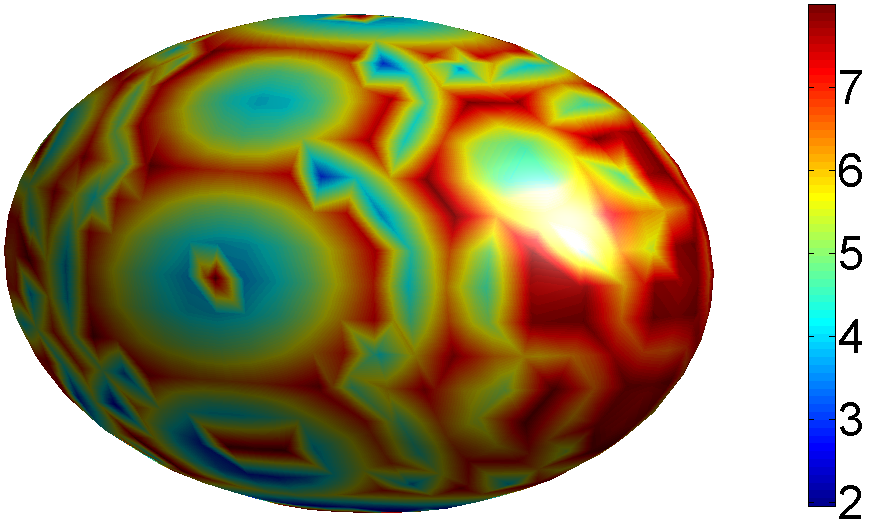} &
\includegraphics[width=0.51\textwidth]{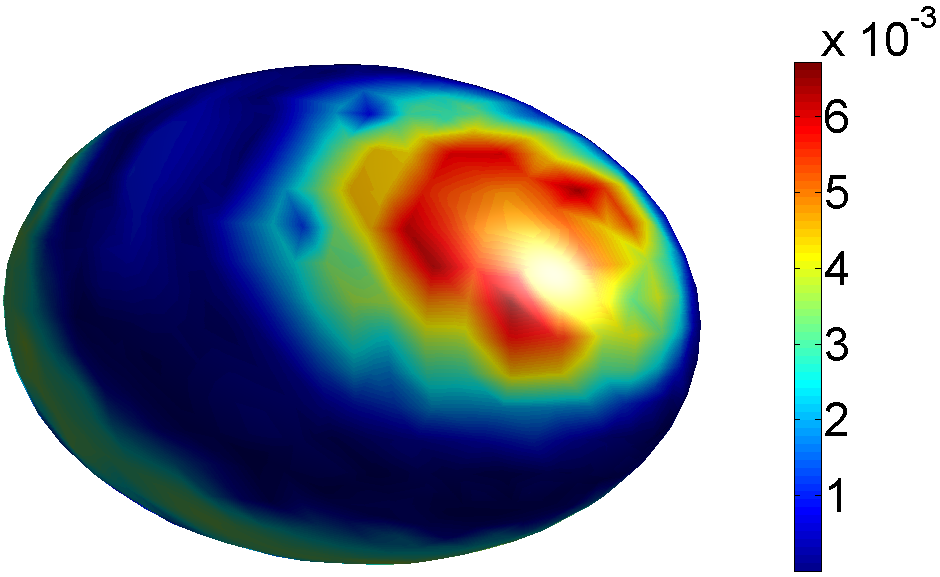}
\end{tabular}
\end{center}
\caption{The computed solution on a  $40 \times 40 \times 40$  mesh (Left chart) and the $L_{\infty}$  error (Right chart) for Case 8(b).}
\label{h2ellipsoid_pic}
\end{figure}

\begin{table}
\caption{Numerical  errors and convergence orders for an ellipsoid interface problem (Case 8(b)).}
\label{h2ellipsoid_tab}
\begin{center}
\begin{tabular}{c||cc|cc}
\hline
$n_x \times n_y \times n_z$ &$L_{\infty}(u)$ &Order   &$L_2 (u )$ &Order     \\\hline
20 $\times$ 20 $\times$ 20   & 3.2934e-2  &       &1.1797e-2   &            \\
40 $\times$ 40 $\times$ 40   & 7.5420e-3  & 2.13  &2.9620e-3   & 1.99        \\
80 $\times$ 80 $\times$ 80   & 1.9183e-3  & 1.98  &7.5382e-4   & 1.97         \\
160 $\times$ 160 $\times$ 160   & 4.7448e-4  & 2.02  &1.8701e-4   & 2.01         \\
\hline
\end{tabular}
\end{center}
\end{table}

In both cases, the second order derivatives of $u_-$ blow up at the origin, making the solution only $H^2$ continuous in $\Omega^-$. Since the differential form, Eq.  (\ref{theoryeq1}), is of second order for the diffusion term, mathematically,  Eq.  (\ref{theoryeq1}) is invalid even though the underlying conservation law    is still valid. However, it is noticed that the equivalent term in the integral form, Eq.  (\ref{integralform1}), involves only the first order derivative. This reduction of the derivative order is important in dealing with solution which changes so rapidly in space that the spatial derivative does not exist. For these reasons, the finite volume method is preferred over the finite difference method in solving problems whose solution has less regularity and exhibits local discontinuities.

In Case 8(a), we choose constant diffusive coefficients $\beta^+$ and $\beta^-$, whereas in Case 8(b), $\beta^+$ and $\beta^-$ are both position dependent. The numerical solutions and the $L_{\infty}$ errors of Case 8(a) and Case 8(b) are depicted in  Figure \ref{h2sphere_pic} and Figure \ref{h2ellipsoid_pic}, respectively. From the numerical results collected in Table \ref{h2sphere_tab} and Table \ref{h2ellipsoid_tab}, it is seen that the second  order accuracy in both cases are essentially obtained.


\section{Conclusion}

In this work, we introduce the finite volume formulation of matched interface and boundary method (MIB-FVM) for solving elliptic interface problems arising from modeling material interface in practice. Much effort has been taken to develop advanced numerical schemes for such  problems in the past few decades. However, challenges still remain in this field. One of the challenge concerns the development of methods with higher order accuracy. Another challenge is to develop methods for dealing with complex interface geometries. The matched interface and boundary (MIB) method has been proved to be able to deal with these challenges. However, since the MIB method is based on the collocation formulation, it does not work well for solutions with low regularity. Based on the integral form rather than the differential form, the finite volume method (FVM) is known for being able to deal with problems with low regularity solutions and better conserve mass and flux. Second order finite volume methods have also been formulated in literature for elliptic interface problems. Motivated by these successes, we propose the present MIB-FVM to take the advantages of MIB and FVM.

To reduce the computational cost  of mesh generation, we utilize  the Cartesian mesh, although the proposed MIB-FVM can be realized on irregular meshes as well. We also employ the vertex-centered FVM, for which the cubic control volumes are associated with the grid point at the center. The control volumes are categorized into regular and irregular types, where special treatments of the MIB formulation are needed for the irregular ones in order to maintain the designed  accuracy. We study the proposed MIB-FVM by a number  of classical test cases, including both 2D cases such as 6-petal flower and jigsaw-puzzle like shape, and 3D cases such as sphere, ellipsoid, standard cylinder, flower-based cylinder and torus. The numerical results all demonstrate the second order accuracy of the proposed method.  The ability of   handling interface with complex geometries and solutions with low regularity ($H^2$ continuous) indicates that the proposed MIB-FVM  combines the advantage of both  MIB and FVM in meeting  numerical challenges. Future work will be done to improve the presnt method in dealing with even more complex interface geometries, for example, interface with geometric  singularities which are commonly seen in practical applications such as electrostatic analysis of proteins.

\section*{Acknowledgments}

This work was supported in part by NSF grants IIS-1302285  and DMS-1160352, NIH Grant R01GM-090208, and
MSU Center for Mathematical Molecular Biosciences Initiative.

\small


\end{document}